\documentclass[notitlepage]{amsart}

\usepackage{amsthm, amssymb, latexsym, tabls, hhline}
\usepackage[vcentermath,enableskew]{youngtab}

\input xypic
\xyoption{all}
\xyoption{arc}

\def\longmapright#1{\hspace{0.3em}\smash{
     \mathop{\longrightarrow}\limits^{#1}}\hspace{0.3em}}

\newtheorem{prop}{Proposition}[section]

\newtheorem{thm}[prop]{Theorem}
\newtheorem{cor}[prop]{Corollary}
\newtheorem{lem}[prop]{Lemma}

\def\equat{\refstepcounter{prop}$$~}
\def\endequat{\leqno{\boldsymbol{(\arabic{section}.\arabic{prop})}}~$$}

\def\exemple#1{{\refstepcounter{prop}\label{#1}\noindent\it Example
\arabic{section}.\arabic{prop} - }}

\def\remarque#1{{\refstepcounter{prop}\label{#1}\noindent\it Remark
\arabic{section}.\arabic{prop} - }}

\def\cba{{\bar{c}}}
\def\iba{{\bar{i}}}

\def\kba{{\bar{k}}}
\def\lba{{\bar{l}}}
\def\mba{{\bar{m}}}
\def\nba{{\bar{n}}}

\def\Cb{{\bs C}}
\def\Pb{{\bs P}}
\def\Qb{{\bs Q}}
\def\infspe{\hspace{0.1em}\mathop{\preccurlyeq}\nolimits\hspace{0.1em}}
\mathchardef\inferieur="321E

\def\ex{\noindent{\it Example. }} % example
\def\rem{\noindent{\it Remark. }} % remarque

\def\getsB{{\buildrel B\over\longleftarrow}}
\def\getsR{{\buildrel R\over\longleftarrow}}

\def\CC{{\mathcal{C}}}
\def\DC{{\mathcal{D}}}
\def\EC{{\mathcal{E}}}

\def\PC{{\mathcal{P}}}
\def\QC{{\mathcal{Q}}}

\def\SC{{\mathcal{A}}}

\def\UC{{\mathcal{U}}}

\def\fonctio#1#2#3#4{\begin{array}{ccc}
{#1} & \longto & {#2} \\
{#3} & \longmapsto & {#4}
\end{array}}

\def\RM{{\mathbb{R}}}
\def\NM{{\mathbb{N}}}
\def\QM{{\mathbb{Q}}}
\def\ZM{{\mathbb{Z}}}

\def\UC{{\mathcal{U}}}
\def\UD{{\mathcal{D}}}

\def\itemth#1{\item[$({\mathrm{#1}})$]}
\def\e{\varepsilon}
\def\t{\tau}

\def\SBT{{\mathcal{SBT}}}

\newcommand{\bs}[1]{{\boldsymbol{#1}}} %bolsymbol
\newcommand{\bbar}[1]{{\overline{#1}}}

\newcommand{\entier}[1]{{[1,#1]}} % entier de 1 a n
\newcommand{\Dentier}[2]{{[#1 ,  #2 ]}} % entier de p a n p<n

\newcommand{\set}[1]{{\left\{ #1 \right\}}}  %ensemble

\def\implies{{\,\Rightarrow\,}} %implies
\def\st{{\,\arrowvert\, }}  %such that

\def\sgn{{\textrm{sign}\,}} %signe d'une lettre, nombre

\def\Sn{{\mathfrak S}} % groupe symetrique
\newcommand{\scalar}[2]{\langle\,#1,#2\,\rangle}
\def\class{\mathcal{CL}_{\mathbb Q}} % algebre des caracteres sur Q
\def\Ind{\mathop{\mathrm{Ind}}\nolimits}
\def\Res{\mathop{\mathrm{Res}}\nolimits}
\def\Irr{\mathop{\mathrm{Irr}}\nolimits}% ensemble des caracteres irreductibles
\def\Ker{\mathop{\mathrm{Ker}}\nolimits}
\def\Rad{\mathop{\mathrm{Rad}}\nolimits}
\def\Cop{\mathop{\mathrm{Cop}}\nolimits}
\def\aug{\mathop{\mathrm{aug}}\nolimits}
\def\rank{\mathop{\mathrm{rank}}\nolimits}
\def\slex{_{\textrm{sl}}} % special lexicographical order
\def\lex{{\mathrm{lex}}} % lexicographical order
\def\VDash{{\, \scriptstyle{ |\models}\,}}  % symbol pour C "signed composition" de n

\def\sh{\mathop{\mathrm{sh}}\nolimits} %shape
\def\cop{\mathcal Q} %coplactic space
\def\Comp{\mathop{\mathrm{Comp}}\nolimits}
\def\Bip{\mathop{\mathrm{Bip}}\nolimits}

\def\SP{{\mathcal{SP}}} % signed permutation

\def\SPprod{{\boldsymbol\ast}} % product on signed permutation
\def\SPcoprod{{\Delta}} % coproduct on signed permutation
\def\Hclass{{\mathcal{CHAR}}}  % Hopf algebra of class functions
\def\CLprod{{\bullet}} % product on class functions
\def\sts{{\textnormal{sts}\,}} %standard signed permutation
\def\alph{{\textnormal{alph}\,}}% set of absolute value of letters

\def\a{\alpha}
\def\b{\beta}
\def\g{\gamma}

\def\d{\delta}
\def\D{\Delta}
\def\e{\varepsilon}
 
\def\l{\lambda}

\def\m{\mu}

\def\s{\sigma}
\def\th{\theta}

\def\t{\tau}
\def\x{\xi}

\def\lamb{{\boldsymbol{\lambda}}}

\def\iff{\Leftrightarrow}

\def\to{\rightarrow}
\def\longto{\longrightarrow}

\def\lexp#1#2{\kern\scriptspace\vphantom{#2}^{#1}\kern-\scriptspace#2}
\def\finl{~$\scriptstyle{\square}$}
\def\SS{\scriptstyle}
\def\DS{\displaystyle}
\def\lamh{{\hat{\l}}}

\def\SG{{\mathfrak S}}
\def\TG{{\mathfrak T}}

\def\eval{{\mathop{\mathrm{ev}}}}

\def\vertical{$\vphantom{\frac{\DS{A}}{\DS{A}}}$}
\def\unbar{{\bar{1}}}
\def\deuxbar{{\bar{2}}}

\def\ch{\textnormal{ch}} %caracteristique de Frobenius

\newcommand{\tab}[2]{{\textnormal{Tab}\left(#1 , #2\right)}}
\newcommand{\bitab}[2]{{\textnormal{Bitab}(#1 , #2)}}

\newcommand\ten{10}
\newcommand\eleven{11}
\newcommand\twelve{12}
\newcommand\thirteen{13}
\newcommand\fourteen{14}
\newcommand\fifteen{15}

\begin{document}

\title{Generalized  descent algebra\\
and construction of  irreducible characters of  hyperoctahedral groups}

\author{C\'edric Bonnaf\'e and Christophe Hohlweg}

\maketitle

%\begin{abstract}
%\end{abstract}

\markboth{}{}

%\tableofcontents

%%
%%
%% Introduction
%%
%%

\section{Introduction}

\medskip

Let $(W,S)$ be a finite Coxeter system and let $\ell : W \to \NM$
denote the length function. If $I\subset S$,   $W_I=<I>$ is the
standard parabolic subgroup generated by $I$ and  $X_I=\{w \in
W~|~ \forall~s \in I,~\ell(ws)>\ell(w)\}$ is a cross-section of
$W/W_I$. Write $x_I = \sum_{w \in X_I} w \in \ZM W$, then
$\Sigma(W)=\oplus_{I \subset S} \ZM x_I$ is a subalgebra of $\ZM
W$ and the $\ZM$-linear map $\th : \Sigma(W) \to \ZM \Irr W$, $x_I
\mapsto \Ind_{W_I}^W 1$ is a morphism of algebras:
 $\Sigma(W)$ is called the {\it descent algebra} or
the {\it Solomon algebra} of $W$ \cite{solomon}.
However, the morphism $\th$ is surjective if and only if
$W$ is a product of symmetric groups.

The aim of this paper is to construct, whenever $W$ is of type $C$,
a subalgebra $\Sigma'(W)$ of $\ZM W$ containing $\Sigma(W)$ and a
surjective morphism of algebras $\th' : \Sigma'(W) \to \ZM \Irr W$
build similarly as $\Sigma(W)$ by starting with a bigger
generating set. More precisely, let $(W_n,S_n)$
denote a Coxeter system of type $C_n$ and write
$S_n=\{t,s_1,\dots,s_{n-1}\}$ where the Dynkin diagram of $(W_n,S_n)$ is
$$
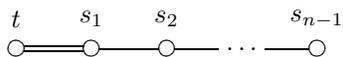
\xy (0,0) *++={} *\frm{o} ; (10,0) *++={} *\frm{o} **@{=}; (20,0)
*++={} *\frm{o} **@{-}; (30,0) *++={\dots}  **@{-}; (40,0) *++={}
*\frm{o} **@{-}; (0,4) *+={t}; %(5,2) *+{4};
(10,4) *+={s_1};
(20,4) *+={s_2};
(40,4) *+={s_{n-1}};% (10,0) *++={>};
\endxy
$$
Let $t_1=t$ and $t_i = s_{i-1} t_{i-1} s_{i-1}$ ($2 \le i \le n$)
and  $S_n'=S_n \cup \{t_1,\dots,t_n\}$. Let $\PC_0(S_n')$ 
denote the set of subsets $I$ of $S_n'$ such that $I=<I> \cap S_n'$. 
If $I \in \PC_0(S_n')$, let $W_I$, $X_I$ and $x_I$ be defined 
as before. Then:

\bigskip

\noindent{\bf Theorem.}
{\it $\Sigma'(W_n)=\oplus_{I \in \PC_0(S_n')} \ZM x_I$
is a subalgebra of $\ZM W_n$ and the $\ZM$-linear map
$\th_n : \Sigma'(W_n) \to \ZM \Irr W_n$,
$x_I \mapsto \Ind_{W_I}^W 1$ is a surjective morphism of algebras.
Moreover, $\Ker \th_n = \sum_{I \equiv I'} \ZM(x_I-x_{I'})$
and $\QM \otimes_\ZM \Ker \th_n$ is the radical of
the $\QM$-algebra $\QM \otimes_\ZM \Sigma'(W_n)$.}

\bigskip

In this theorem, the notation $I \equiv I'$ means that
there exists $w \in W_n$ such that $I' = \lexp{w}{I}$, that is, $W_I$ and $W_{I'}$
 are conjugated.
This theorem is stated and proved in \S\ref{main result}. Note
that it is slightly differently formulated: in fact, it turns out
that there is a natural bijection between signed compositions of
$n$ and $\PC_0(S_n')$ (see Lemma \ref{reconnaissance}). So,
everything in the text is indexed by signed compositions instead
of $\PC_0(S_n')$. It must also be noticed that, by opposition with
the classical case, the multiplication $x_I x_J$ may involve
negative coefficients. 
Using another basis, we show that $\Sigma'(W_n)$ is precisely the 
generalized descent algebra discovered by Mantaci and Reutenauer 
\cite{mantaci-reutenauer}.

Using this theorem and the Robinson-Schensted correspondence for
type $C$ constructed by Stanley \cite{stanley} and a Knuth version
of it given in \cite{bonnafe-iancu}, we obtain an analog of
J\"ollenbeck's result (on the construction of characters of the
symmetric group \cite{jollenbeck}) using an extension
$\tilde{\th}_n : \QC_n \to \ZM\Irr W_n$ of $\th_n$ to the
coplactic space $\QC_n$ (see Theorem \ref{theoreme theta}). The
coplactic space refer to J\"ollenbeck's construction revised in
\cite{bless}.

Now,  let $\SP=\oplus_{n \ge 0} \ZM W_n$, $\Sigma'=\oplus_{n \ge
0} \Sigma'(W_n)$ and $\QC=\oplus_{n \ge 0} \QC_n$. Let
$\th=\oplus_{n \ge 0} \th_n$ and $\tilde{\th}=\oplus_{n \ge 0}
\tilde{\th}_n$. Aguiar and Mahajan have proved that $\SP$ is
naturally a Hopf algebra and that $\Sigma'$ is a Hopf subalgebra
\cite{aguiar-mahajan}. We prove here that $\QC$ is also a Hopf
subalgebra of $\SP$ (containing $\Sigma'$) and that $\th$ and
$\tilde{\th}$ are surjective morphisms of Hopf algebras (see
Theorem \ref{hopf_coplac}). This generalizes similar results in
symmetric groups (\cite{poirier-reutenauer} and \cite{bless}),
which are parts of  combinatorial tools used within the framework
of the representation theory of type $A$ (see for instance
\cite{thibon}).

In the last section of this paper, we give some explicit
computations in $\Sigma'(W_2)$ (characters, complete set
of orthogonal primitive idempotents, Cartan matrix of $\Sigma'(W_2)$...). \\

In the Appendix, P.~Baumann and the second author link the above 
construction  with the Specht construction and symmetric functions 
(see \cite{macdonald}).

\bigskip

\noindent{\it Remark.}
It seems interesting to try to construct
a subalgebra $\Sigma'(W)$ of $\ZM W$ containing $\Sigma(W)$
and a morphism $\th' : \Sigma'(W) \to \ZM \Irr W$
for arbitrary Coxeter group $W$.
But it is impossible to do so in a same way as we did for
type $C$ (by extending the generating set). Computations
using {\tt CHEVIE} programs show us that it is impossible
to do so in type $D_4$ and that the reasonable choices
in $F_4$ fail (we do not obtain a subalgebra!).
However, it is possible to do something similar for type $G_2$.
More precisely, let $(W,S)$ be of type $G_2$. Write
$S=\{s,t\}$ and let $S'=\{s,t,sts,tstst\}$ and repeat the
procedure described above to obtain a sub-$\ZM$-module
$\Sigma'(W)$ of $\ZM W$ and a morphism $\th' : \Sigma'(W) \to \ZM\Irr W$.
Then the theorem stated in this introduction also holds in this case.
We have $\mathop{\mathrm{rank}}_\ZM \Sigma'(W)=8$ and
$\mathop{\mathrm{rank}}_\ZM \Ker \th'=2$.

%
%
% Some reflection subgroups of hyperoctahedral groups
%

%\section{Generalized  descent algebras}\label{gene descent algebra}

\section{Some reflection subgroups of hyperoctahedral groups}

In this article,
we denote $[m,n]=\{i \in \ZM~|~m \le i \le n\}=\{m,m+1,\dots,n-1,n\}$, for all
$m\leq n\in\ZM$, and $\sgn (i)\in\set{\pm 1}$ the
sign of $i\in\ZM\setminus\{0\}$.  If $E$ is a set, we denote by
%$\PC(E)$ the set of subsets of $E$ and by
$\SG(E)$ the group of permutations on the set $E$.
If $m \in \ZM$, we often denote by $\bbar{m}$
the integer $-m$.

%
% hyperoctahedral group
%

\subsection{The hyperoctahedral group} We begin by making clear some notations
and definitions concerning the {\it hyperoctahedral group} $W_n$.
%Let $(W_n,S_n)$ be a Coxeter system of type $C_n$,  where
% $S_n=\set{t, s_1 , \dots , s_{n-1}}$ is such that t and
%$t_{i+1}=s_i t_i s_i$ for any  $1\leq i\leq n-1$.he Coxeter graph of $W_n$ is
%$$
%\xy (0,0) *++={} *\frm{o} ; (10,0) *++={} *\frm{o} **@{=}; (20,0)
%*++={} *\frm{o} **@{-}; (30,0) *++={\dots}  **@{-}; (40,0) *++={}
%*\frm{o} **@{-}; (0,4) *+={t}; %(5,2) *+{4};
%(10,4) *+={s_1};
%(20,4) *+={s_2};
%(40,4) *+={s_{n-1}};% (10,0) *++={>};
%\endxy
%$$
Denote
%$\ell(w)$ the length of $w\in W_n$ as a word in the
%elements of $S_n$, and
$1_n$  the identity of   $W_n$ (or $1$ if no confusion is
possible). We denote by $\ell_t(w)$ the number of occurrences of
$t$ in a reduced decomposition of $w$ and we define
$\ell_s(w)=\ell(w)-\ell_t(w)$.

It is well-known that $W_n$ acts on the set $I_n=[1,n] \cup [\bar n,\bar 1]$
by permutations as follows: $t=(\bar 1,\ 1)$ and $s_i = (\bbar{i+1}, \bbar i)(i, i+1)$
for any $i\in [1,n-1]$. Through this action, we have
$$
W_n=\{w \in \SG(I_n)~|~\forall~i \in I_n,~w(\bbar{\, i\,})=\bbar{w(i)}\} .
$$
We often represent $w\in W_n$ as the word $w(1)w(2)\dots w(n)$ in examples.

The subgroup $W_\nba=\{w \in W_n~|~w([1,n]) = [1,n]\}$ of $W_n$ is
naturally identified with $\SG_n$, the symmetric group of degree
$n$, by restriction of its elements to $[1,n]$. Note that $W_\nba$
is generated, as a reflection subgroup of $W_n$, by $S_\nba =
\{s_1,\dots, s_{n-1}\}$.

A {\it standard parabolic subgroup} of $W_n$ is a subgroup
generated by a subset of $S_n$ (a {\it parabolic subgroup} of
$W_n$ is a subgroup conjugate to some standard parabolic subgroup).
Note that $(W_\nba,S_\nba)$ is a Coxeter group, which is a
standard parabolic subgroup of $W_n$. If $m \le n$, then $S_m$ is
naturally identified with a subset of $S_n$ and $W_m$ will be
identified with the standard parabolic subgroup of $W_n$ generated
by $S_m$.\\

Now, we set $T_n=\{t_1,\dots,t_n\}$, with $t_i$ as in
Introduction.
 As a permutation of $I_n$,
note that $t_i=(i,\bar{i})$, then the
reflection subgroup $\TG_n$ generated by $T_n$ is naturally
identified with $(\ZM/2\ZM)^n$. Therefore $W_n=W_\nba \ltimes
\TG_n$ is just the wreath product of $\SG_n$ by  $\ZM/2\ZM$. If $w
\in W_n$, we denote by $(w_S,w_T)$ the unique pair in $\SG_n
\times \TG_n$ such that $w=w_S w_T$. Note that
$\ell_t(w)=\ell_t(w_T)$. In this article, we will  consider
reflection subgroups generated by subsets of $S_n'=S_n \cup T_n$.

%
% root system
%

\subsection{Root system.}
Before studying the reflection subgroups generated by subsets of $S'_n$,
let us  recall some basic facts about Weyl groups of type $C$ (see \cite{bourbaki}).
 Let us endow $\RM^n$ with its canonical euclidean structure.
Let $(e_1,\dots,e_n)$  denote the canonical basis of $\RM^n$: this is an orthonormal
basis. If $\alpha \in \RM^n\setminus\{0\}$, we denote by $s_\alpha$
the orthogonal reflection such that $s_\a(\a)=-\a$.   Let
$$\Phi_n^+=\{2 e_i~|~1 \le i \le n\} \cup \{e_j + \nu e_i~|~
\nu \in \{1,-1\}~\text{and}~1 \le i < j \le n \},$$
$\Phi_n^-=-\Phi_n^+$ and $\Phi_n = \Phi^+_n \cup \Phi_n^-$. Then $\Phi_n$ is a
root system of type $C_n$ and $\Phi_n^+$ is a positive root system of $\Phi_n$.
By sending $t$ to $s_{2e_1}$ and $s_i$ to $s_{e_{i+1}-e_i}$
(for $1 \le i \le n-1$), we will identify $W_n$ with the Coxeter
group of $\Phi_n$. Then
$$
\D_n=\{2e_1,e_2-e_1,e_3-e_2,\dots,e_n-e_{n-1}\}
$$
is the basis of $\Phi_n$ contained in $\Phi_n^+$ and
the subset $S_n$ of $W_n$ is naturally
identified with the set of simple reflections
$\{s_\a~|~\a \in \D_n\}$. Therefore, for any $w\in W_n$ we have
$$\ell(w)=|\Phi_n^+\cap w^{-1}(\Phi_n^-)|;$$
and $\ell(ws_\alpha)<\ell(w)$ if and only if
$w(\alpha)\in\Phi^-$, for all  $\alpha\in\Phi^+$.\\

\remarque{superieur} Let $w \in W_n$ and let $\a \in \Phi_n^+$. Then
$\ell(ws_\alpha)<\ell(w)$ if and only if
$w(\alpha)\in\Phi_n^-$.
Therefore, if $i \in [1,n-1]$, then
$$
\ell(ws_i)<\ell(w) \iff w(i)>w(i+1),
$$
and, if $j \in [1,n]$, then
$$
\ell(wt_j)<\ell(w)\iff w(j)<0.
$$
Therefore, we deduce from the strong exchange condition (see
\cite[\S5.8]{humphreys})
 \equat\label{t longueur}
\ell_t(w)=|\{i \in [1,n]~|~w(i) < 0 \}|.
\endequat

%In this description, a subgroup of $W_n$ is parabolic if and only if it
% is the stabilizer of some element of $\RM^n$.

\medskip

\subsection{Some closed subsystems of $\Phi_n$}\label{reflection subgroups}

Consider the subsets $\set{s_1,t_1}$ and $\set{s_1,t_2}$ of $S'_n$
($n\geq 2)$. It is readily seen that these two sets of reflections
generate the same reflection subgroup of $W_n$. This lead us to
find a parametrization of subgroups generated by a subset of $S'_n$.

A {\it signed composition} is a
sequence $C = (c_1, \dots, c_r)$  of non-zero elements of $\ZM$. The number $r$ is
called the {\it length} of $C$. We set $|C| =\sum_{i=1}^r | c_i |$. If $|C|=n$,
we say that $C$ is a {\it signed composition of $n$} and we write $C \VDash n$.
We also define
$C^+=(|c_1|,\dots,|c_r|) \VDash n$, $C^-=-C^+$ and  $\bbar C = -C$. We denote
by $\Comp(n)$ the set of signed compositions of $n$. In
particular, any composition is a signed composition (any part is
positive). Note that \equat\label{nombre composition}
|\Comp(n)|=2.3^{n-1}.
\endequat

Now, to each $C=(c_1, \dots, c_r) \VDash n$, we associate a reflection subgroup of
$W_n$ as follows:  for $1\leq i\leq r$, set
$$I_C^{(i)}=\begin{cases}
            I_{C,+}^{(i)} & \text{if } c_i < 0,\\
        I_{C,+}^{(i)} \cup -I_{C,+}^{(i)} & \text{if } c_i > 0,
        \end{cases}$$
where $I_{C,+}^{(i)}=\big[|c_1|+\dots+|c_{i-1}|+1,|c_1|+\dots+|c_i|\big]$. Then
$$W_C=\{w \in W_n ~|~ \forall~1 \le i \le r,~w(I_C^{(i)})=I_C^{(i)}\}$$
is a reflection subgroup generated by
\begin{eqnarray*}
S_C& =&  \set{s_p \in S_\nba ~|~
|c_1| + \dots+|c_{i-1}|+1\leq p \leq |c_1| +\dots+ |c_i|  -1}\\
&&\cup\set{t_{|c_1| + \dots+|c_{j-1}|+1} \in T_n \st
c_j >0}\quad\subset S'_n
\end{eqnarray*}
Therefore, $W_C \simeq W_{c_1} \times \dots \times W_{c_r}$: we
denote by $(w_1,\dots,w_r) \mapsto w_1 \times \dots \times w_r$
the natural isomorphism $W_{c_1} \times \dots \times W_{c_r} \longmapright{\sim} W_C$.

\medskip

\noindent{\it Example.}
The group $W_{(\bbar 2,  3,\bbar 1,\bbar 3,  1)}
\simeq \Sn_2\times W_3 \times \Sn_1 \times \Sn_3 \times W_1$
is generated, as a reflection subgroup of $W_{10}$, by
$S_{(\bbar 2,  3,\bbar 1,\bbar 3,  1)}=
\set{s_1 }\cup\set{t_3, s_3, s_4}\cup \set{s_7,s_8}\cup\set{t_{10}}\subset S'_{10}$.

\medskip

The signed composition $C$ is said
{\it semi-positive} if $c_i \geq -1$ for every $i \in [1,r]$. Note that
a composition is a semi-positive composition. We say that $C$
is {\it negative} if $c_i < 0$ for every $i \in [1,r]$. We say
that $C$ is {\it parabolic} if $c_i < 0$ for $i \in [2,r]$. Note
that $C$ is parabolic if and only if $W_C$ is a standard parabolic
subgroup.

\medskip

Now, let $S_C'=S_n' \cap W_C$, $\Phi_C=\{\a \in \Phi_n~|~s_\a \in
W_C\}$ and $\Phi_C^+=\Phi_C \cap \Phi_n^+$. Then $W_C$ is the Weyl
group of the closed subsystem $\Phi_C$ of $\Phi_n$. Moreover,
$\Phi_C^+$ is a positive root system of $\Phi_C$ and we denote by
$\D_C$ the basis of $\Phi_C$ contained in $\Phi_C^+$. Note that
$S_C=\{s_\a~|~\a \in \D_C\}$, so $(W_C,S_C)$ is a Coxeter group.

Let $\ell_C : W_C \to \NM$
denote the length function on $W_C$ with respect to $S_C$. Let $w_C$ denote
the longest element of $W_C$ with respect to
$\ell_C$. If $C$ is a composition, we denote by $\s_C$ the longest
element of $\SG_C=W_{\bbar C}$ with respect to $\ell_{\bbar C}$
(which is the restriction of $\ell$ to $\SG_C$). In other words,
$\s_C=w_{\bbar C}$. In particular, $w_n$ (resp. $\s_n$) denotes
the longest element of $W_n$ (resp. $\SG_n$).

Write $T_C=T_n \cap W_C$ and  $\TG_C=\TG_n \cap W_C$, then observe
that
\equat
W_C = W_{C^-} \ltimes \TG_C = \SG_{C^+} \ltimes \TG_C.
\endequat

\medskip

\noindent{\it Remarks.} (1) This class of reflection subgroups contains
the standard
parabolic subgroups, since $S_n\subset S'_n$. But it contains also
some other subgroups which are not parabolic (consider the
subgroup generated by $\set{t_1,t_2}$ as example). In other words,
it may happen that $\D_C \not\subset \D_n$. In fact, $\D_C \subset
\D_n$ if and only if $W_C$ is a standard parabolic subgroup of
$W_n$.

(2) If $W_C$ is not a standard parabolic subgroup of $W_n$, then
$\ell_C$ is not the restriction of $\ell$ to $W_C$.

\bigskip

We close this subsection by an easy characterization of the subsets $S_C'$:

\begin{lem}\label{reconnaissance}
Let $X$ be a subset of $S_n'$. Then the following are equivalent:
\begin{itemize}
\itemth{1} $<X> \cap S_n' = X$.

\itemth{2} $X \cap T_n$ is stable under conjugation by $<X>$.

\itemth{3} $X \cap T_n$ is stable under conjugation by $<X \cap S_\nba>$.

\itemth{4} There exists a signed composition $C$ of $n$ such that $X=S_C'$.
\end{itemize}
\end{lem}

\begin{cor}\label{conjugaison s'}
Let $w \in W_n$ and let $C \VDash n$. If $\lexp{w}{S_C'} \subset S_n'$,
then there exists a (unique) signed composition $D$ such that
$\lexp{w}{S_C'}=S_D'$.
\end{cor}

\begin{proof}
Indeed, $\lexp{w}{S_C'} \cap T_n = \lexp{w}{(S_C' \cap T_n)}$ and
$\lexp{w}{S_C'} \cap S_\nba = \lexp{w}{(S_C' \cap S_\nba)}$.
\end{proof}

%\subsubsection{Corresponding generalized Coxeter classes}

%
% Orbits of closed subsystems of ${\boldsymbol{\Phi_n}}$
%

\subsection{Orbits of closed subsystems of $\Phi_n$}
In this subsection,
we determine when two subgroups $W_C$ and $W_D$
of $W_n$ are conjugated.
A \textit{bipartition} of $n$ is a pair
$\l= (\l^+ , \l^-)$ of partitions such that $|\l|:=|\l^+|+|\l^-|=n$.
We write $\l \Vdash n$ to say that $\l$ is a bipartition of $n$, and
the set of bipartitions of $n$ is denoted by
$\Bip(n)$. It is well-known that the conjugacy classes of $W_n$
are in bijection with $\Bip(n)$ (see \cite{geck,macdonald}). We
define ${\hat{\l}}$ as the signed composition of $n$ obtained by
concatenation of $\l^+$ and $-\l^-$. The map $\Bip(n) \to
\Comp(n)$, $\l \mapsto \hat{\l}$ is injective.

Now, let $C$ be a signed composition of $n$. We define
$\lamb(C)=(\l^+,\l^-)$ as the bipartition of $n$ such that $\l^+$
(resp. $\l^-$) is obtained from $C$ by reordering in decreasing
order the positive parts of $C$ (resp. the absolute value of the
negative parts of $C$). One can easily check that the map
$$\lamb : \Comp(n) \longto \Bip(n)$$
is surjective (indeed, if $\l \in \Bip(n)$, then
$\lamb(\hat{\l})=\l$) and that the following proposition holds:

\begin{prop}\label{conjugaison}
Let $C,D\VDash n$, then $W_C$ and $W_D$ are conjugate in $W_n$ if
and only if $\lamb(C)=\lamb(D)$. If $\Psi$ is a closed subsystem
of $\Phi_n$, then there exists a unique bipartition $\l$ of $n$
and some $w \in W_n$ such that $\Psi=w(\Phi_{\hat{\l}})$.
\end{prop}

Let $C, D\VDash n$, then we write $C \subset D$ if $W_C \subset
W_D$. Moreover, $C$, $C' \subset D$ and if $W_C$ and $W_{C'}$ are
conjugate under $W_D$, then we write $C \equiv_D C'$.

%
% Cosets
%

\subsection{Distinguished coset representatives} Let $C \VDash
n$, then
$$
X_C=\{x \in W_n~|~\forall~w \in W_C,~\ell(xw) \ge \ell(x)\}
$$
is  a  distinguished set of {\it minimal coset representatives}
for $W_n/W_C$ (see proposition below). It is readily seen that
\begin{eqnarray*}
X_C&=&\{w \in W_n~|~w(\Phi_C^+) \subset \Phi_n^+\}\\
&=&\{w \in W_n~|~\forall ~\a \in \D_C,~w(\a) \in \Phi_n^+\}.
\end{eqnarray*}
Finally
$$
X_C= \{w \in W_n~|~\forall~r \in S_C,~\ell(wr)>\ell(w)\}.
$$

We need a relative notion: if $D\VDash n$ such that $C\subset D$,
the set $X_C^{D}=X_C \cap W_{D}$ is a distinguished set of
minimal coset representatives for $W_D/W_C$. If $D=(n)$
we write $X_C^n$ instead of $X_C^{(n)}$.

\begin{prop}\label{simple classe}
Let $C\VDash n$, then:
\begin{itemize}
\itemth{a} The map $X_C \times W_C \to W_n$, $(x,w) \mapsto xw$ is
bijective.

\itemth{b} If $C \subset D$, then the map $X_D \times X_C^D \to
X_C$, $(x,y) \mapsto xy$ is bijective.

\itemth{c} If $x \in X_C$ and $w \in W_C$, then $\ell_t(xwx^{-1})
\ge \ell_t(w)$. Consequently, $\SG_n \cap \lexp{x}{W_C} = \SG_n
\cap \lexp{x}{\SG_{C^+}}$.
\end{itemize}
\end{prop}

\begin{proof}
(a) is stated, in a general case, in \cite[Proposition 3.1]{fleisch}. (b) follows easily from (a). Let us
now prove (c). Let $x \in X_C$ and $w \in W_C$. Let $I=\{i \in
I_n~|~w(i) < 0\}$ and $J=\{i \in I_n~|~xwx^{-1}(i) < 0 \}$, then
 $\ell_t(w)=|I|$ and $\ell_t(xwx^{-1})=|J|$, by \ref{t longueur}.
Now let $i \in I$, then $t_i \in W_C$, so $\ell_t(xt_i) >
\ell_t(x)$. In other words, $x(i) > 0$. Now, we have
$xwx^{-1}(x(i))=xw(i)$. But, $w(i) < 0$ and $t_{-w(i)}=wt_i w^{-1}
\in W_C$. Therefore, $x(-w(i))=-xw(i) > 0$. This shows that $x(i)
\in J$. So, the map $I \to J$, $i \mapsto x(i)$ is well-defined
and clearly injective, implying $|I| \le |J|$ as desired.

The last assertion of this proposition follows easily from this
inequality and from the fact that $\Sn_{C^+}=\{w \in
W_C~|~\ell_t(w)=0\}$.
\end{proof}

\begin{prop}\label{conjugaison x}
Let $C \VDash n$ and $x \in X_C$ be such that $\lexp{x}{S_C'} \subset
S_n'$. Let $D$ be the unique signed composition of $n$ such that
$\lexp{x}{S_C'} = S_D'$ (see Corollary \ref{conjugaison s'}).
Then $X_C=X_D x$.
\end{prop}

\begin{proof}
By symmetry, it is sufficient to prove that, if $w \in X_D$,
then $wx \in X_C$. Let $\a \in \Phi_C^+$. Then, since $x \in X_C$,
we have $x(\a) \in \Phi_n^+ \cap \lexp{x}{\Phi_C}=\Phi_D^+$. So
$w(x(\a)) \in \Phi_n^+$ since $w \in X_D$. So $wx \in X_C$.
\end{proof}

\medskip

\subsection{Maximal element in $X_C$}
It turns out that, for every signed composition $C$
of $n$, $X_C$ contains a unique element of maximal length (see
Proposition \ref{X maximal}). First, note the following two examples:

(1)   if $C$ is parabolic, it is well-known that $\ell_C$ is the restriction
of $\ell$
and that, for all $(x,w)\in X_C \times W_C$, we have
%\equat\label{addition longueur parabolique}
$$
\ell(xw) = \ell(x) +\ell(w)
$$
%\endequat
In particular, $w_n w_C$ is the longest element of $X_C$ (see
\cite{geck});

 (2) let $C$ be a composition of $n$, then $W_C$ is not in general
a standard parabolic subgroup of $W_n$. However, since $W_C$
contains $\TG_n$, $X_C$ is contained in $\SG_n$. This shows that
$$X_C=X_{\bbar C}^\nba = X_{\bbar C} \cap \SG_n.$$
In particular, $X_C$ contains a unique element of maximal length:
this is $\s_n \s_{C}$;

\bigskip

Now,  let $k$ and $l$ be two non-zero natural numbers such that
$k+l=n$. Then $W_{k,l}$ is not a parabolic subgroup of $W_n$.
However, $W_{k,\lba}$ is a standard parabolic subgroup of $W_n$
and $W_{k,\lba} \subset W_{k,l}$. So $X_{k,l} \subset X_{k,\lba}$.
So, if $x \in X_{k,l}$ and $w \in W_{k,\lba}$, then
\equat\label{addition longueur particulier}
\ell(xw)=\ell(x)+\ell(w).
\endequat
This applies for instance if $w \in W_k \subset W_{k,\lba}$.

Then,  we need to introduce a decomposition of $X_C$ using
Proposition \ref{simple classe} (b). Write $C=(c_1,\dots,c_r)\VDash n$. We
set
$$X_{C,i}=
X_{(|c_1|+\dots+|c_{i-1}|,c_i,\dots,c_r)}^{(|c_1|+\dots+|c_i|,c_{i+1},\dots,c_r)}.$$
Then the map
$$\fonctio{X_{C,r} \times \dots \times X_{C,2} \times X_{C,1}}{X_C}{
(x_r,\dots,x_2,x_1)}{x_r\dots x_2 x_1}$$ is bijective by
Proposition \ref{simple classe} (b). Moreover, by \ref{addition
longueur particulier}, we have \equat\label{addition longueur X}
\ell(x_r\dots x_2 x_1)=\ell(x_r) + \dots + \ell(x_2) + \ell(x_1)
\endequat
for every $(x_r,\dots,x_2,x_1) \in X_{C,r} \times \dots \times
X_{C,2} \times X_{C,1}$. For every $i \in [1,r]$, $X_{C,i}$
contains a unique element of maximal length (see (1)-(2) above). Let us denote it by
$\eta_{C,i}$. We set:
$$\eta_C=\eta_{C,r} \dots \eta_{C,2} \eta_{C,1}.$$
Then, by \ref{addition longueur X}, we have

\begin{prop}\label{X maximal}
Let $C\VDash n$, then $\eta_C$ is the unique element of $X_C$ of
maximal length.
\end{prop}

%
% Double cosets representatives
%

\subsection{Double cosets representatives}
If $C$ and $D$ are two signed compositions of $n$, we set
$$X_{CD}=X_C^{-1} \cap X_D.$$

\begin{prop}\label{double classe}
Let $C$ and $D$ be two signed composition of $n$ and let $d \in
X_{CD}$. Then:
\begin{itemize}
\itemth{a} There exists a unique signed composition $E$ of $n$
such that $S_E'=S_C' \cap \lexp{d}{S_D}'$. It will be denoted by
$C \cap \lexp{d}{D}$ or $\lexp{d}{D} \cap C$. We have $(C \cap
\lexp{d}{D})^-=C^- \cap \lexp{d}{D^-}$.

\itemth{b} $W_C \cap \lexp{d}{W_D}=W_{C \cap \lexp{d}{D}}$ and
$W_C \cap \lexp{d}{S_D'}=S_C' \cap \lexp{d}{W_D} = S_{C \cap
\lexp{d}{D}}'$.

\itemth{c} If $w \in W_{C \cap \lexp{d}{D}}$, then
$\ell_t(w)=\ell_t(d^{-1}w d)$.

\itemth{d} If $w \in W_C d W_D$, then there exists a unique pair
$(x,y) \in X_{C \cap \lexp{d}{D}}^C \times W_D$ such that $w=xdy$.

\itemth{e} Let $(x,y) \in X_{C \cap \lexp{d}{D}}^C \times W_D$,
then $\ell(xdy) \ge \ell(x_S) + \ell_t(x) + \ell(d) + \ell(y_S) +
\ell_t(y)$.

\itemth{f} $d$ is the unique element of $W_C d W_D$ of minimal
length.
\end{itemize}
\end{prop}

\begin{proof}
(a) follows immediately from Lemma \ref{reconnaissance}
(equivalence between (3) and (4)).

%Since $S_{\nba}$ and $T_n$ are contained in two different
%conjugacy classes of $W_n$, we have
%$$S_C' \cap \lexp{d}{S}_D' = (S_{C^-} \cap \lexp{d}{S_{D^-}})
%\coprod (T_C \cap \lexp{d}{T_D}).$$ Since $-C^+$ is parabolic, and
%since $d \in X_{C^-,D^-}$, there exists a unique composition
%$E^+$ of $n$ such that $S_{C^-} \cap \lexp{d}{S}_{D^-} =
%S_{E^-}$. Note that $\SG_{C^+} \cap \lexp{d}{\SG}_{D^+} =
%\SG_{E^+}$.
%
%Let us write $E^+=(e_1^+,\dots,e_r^+)$. If $1 \le k \le r$, we
%define $\e_k \in \{1,-1\}$ as follows: $\nu_k = 1$ if and only if
%there exists $i \in I_{E^+,+}^{(k)}$ such that $t_i \in T_C \cap
%\lexp{d}{T_D}$. Now, let $e_k = \nu_k e_k^+$ and let
%$E=(e_1,\dots,e_r)$. Note that $I_{E^+,+}^{(k)}=I_{E,+}^{(k)} \cap
%\NM$.
%
%Let $k$ be such that $\nu_k=1$. Then there exists $i \in
%I_{E,+}^{(k)}$ such that $t_i \in T_C \cap \lexp{d}{T_D}$. Let now
%$j \in I_{E,+}^{(k)}$. Then there exists $w \in \SG_{E^+}$ such
%that $w t_i w^{-1} = t_j$. But, $w \in \SG_{C^+}$, so $t_j \in
%T_C$ and $w \in \lexp{d}{\SG}_{D^+}$, so $t_j \in \lexp{d}{T_D}$.
%This shows that $\nu_k=-1$ if and only if $t_j \in T_C \cap
%\lexp{d}{T_D}$ for every $j \in I_{E^+}^{(k)}$. In particular,
%$$S_E'=S_C' \cap \lexp{d}{S_D}'.$$

\medskip

(b) It is clear that $W_E \subset W_C \cap \lexp{d}{W_D}$. Let us
show the reverse inclusion. Let $w \in W_C \cap \lexp{d}{W_D}$. We
will show by induction on $\ell_t(w)$ that $w \in W_E$. If
$\ell_t(w)=0$, then we see from Proposition \ref{simple classe}
(d) that $w \in \SG_{C^+} \cap \lexp{d}{\SG_{D^+}}=\SG_{E^+}$ by
definition of $E^+$.

Assume now that $\ell_t(w) > 0$ and that, if $w' \in W_C \cap
\lexp{d}{W_D}$ is such that $\ell_t(w') < \ell_t(w)$, then $w' \in
W_E$. Since $\ell_t(w)>0$, there exists $i \in [1,n]$ such that
$w(i) < 0$. In particular, $t_i \in \TG_C$. By the same argument
as in the proof of Proposition \ref{simple classe} (d), we have
that $t_i \in \lexp{d}{W_D}$. So, $t_i \in T_C \cap \lexp{d}{T_D}
= T_E$. Now, let $w'=wt_i$. Then $t_i \in W_E$, $w' \in W_C \cap
\lexp{d}{W_D}$ and $\ell_t(w')=\ell_t(w)-1$. So, by the induction
hypothesis, $w' \in W_E$, so $w \in W_E$.

The other assertions of (b) follow easily.

\medskip

(c) Let $w=\s_1 \dots \s_l$ be a reduced decomposition of $w$ with
respect to $S_C$. Then $d^{-1} w d =(d^{-1} \s_1 d) \dots (d^{-1}
\s_l d)$. But $d^{-1} \s_i d \in \lexp{d^{-1}}{(S_C' \cap
\lexp{d}{S_D'})}= S_{\lexp{d^{-1}}{C} \cap D}'$, so $\ell_t(d^{-1}
\s_i d)=\ell_t(\s_i)$. Since $\ell_t(w)=\ell_t(\s_1) + \dots +
\ell_t(\s_l)$, we see that $\ell_t(w) \ge \ell_t(d^{-1} w d)$. By
symmetry, we obtain the reverse inequality.

\medskip

(d) Let $w \in W_C d W_D$. Let us write $w=adb$, with $a \in W_C$
and $b \in W_D$. We then write $a=xa'$ with $x \in X_{C \cap
\lexp{d}{D}}^C$ and $a' \in \Sn_{C \cap \lexp{d}{D}}$. Then
$d^{-1} a' d \in W_{\lexp{d^{-1}}{C} \cap D} \subset W_D$. Write
$y=(d^{-1} a' d) b$. Then $(x,y) \in X_{C \cap \lexp{d}{D}}^C
\times W_D$ and $w=xdy$.

Now let $(x',y') \in X_{C \cap \lexp{d}{D}}^C \times W_D$ such
that $w=x'dy'$. Then $x^{\prime -1} x = d(yy^{\prime -1})d^{-1}$.
So $x^{\prime -1} x \in W_{C \cap \lexp{d}{D}}$, that is $x W_C =
x' W_C$. So $x=x'$ and $y=y'$.

\medskip

(e) Let $(x,y) \in X_{C \cap \lexp{d}{D}}^C \times W_D$. We will
show by induction on $\ell_t(x)+\ell_t(y)$ that
$$\ell(xdy) \ge \ell(x_S) + \ell_t(x) + \ell(d) + \ell(y_S) + \ell_t(y).$$
If $\ell_t(x)=\ell_t(y)=0$, then $x \in X^{C^-}_{C^- \cap
\lexp{d}{(D^-)}}$, $y \in \SG_{D^+}$ and $d \in X_{C^-,D^-}$.
So, by \cite[Lemma 2]{bbht}, we have $\ell(xdy) = \ell(x_S)  + \ell(d) +
\ell(y_S)$, as desired.

Now, let us assume that $\ell_t(x)+\ell_t(y) > 0$ and that the
result holds for every pair $(x',y') \in X_{C \cap \lexp{d}{D}}^C
\times W_D$ such that $\ell_t(x')+\ell_t(y') < \ell_t(x) +
\ell_t(y)$. By symmetry, and using (c), we can assume that
$\ell_t(y) > 0$. So there exists $i \in I_n$ such that $y(i) < 0$.
Let $y'=yt_i$. Then $t_i \in T_D$, $\ell(y_S)=\ell(y_S')$,
$\ell_t(y')=\ell_t(y)-1$. Therefore, by induction hypothesis, we
have
$$\ell(xdy') \ge \ell(x_S) + \ell_t(x) + \ell(d) + \ell(y_S) + \ell_t(y)-1.$$
It is now enough to show that $\ell(xdy't_i) > \ell(xdy')$, that
is $xdy'(i) > 0$. Note that $y'(i) > 0$ and that $t_{y'(i)}=y' t_i
y^{\prime -1} \in W_D$. So the result follows from the following
lemma:

\begin{quotation}
\begin{lem}\label{easy}
If $d \in X_{CD}$, if $x \in X_{C \cap \lexp{d}{D}}^D$ and if $j
\in [1,n]$ is such that $t_j \in T_D$, then $xd(j) > 0$.
\end{lem}

\begin{proof}
Since $t_j \in W_D$ and $d \in X_D$, we have $d(j) > 0$. Two cases
may occur. If $t_{d(j)} \in T_C$, then $t_{d(j)}=dt_j d^{-1} \in
T_{C \cap \lexp{d}{D}}$. Therefore, $x(d(j)) > 0$ since $x \in
X_{C \cap \lexp{d}{D}}^C$. If $t_{d(j)} \not\in T_C$, then
$x(d(j)) > 0$ since $x \in W_C=\SG_{C^+} \ltimes \TG_C$.
\end{proof}
\end{quotation}

\medskip

(f) follows immediately from (e).
\end{proof}

\medskip

\remarque{xcd}  Let $C$ and $D$ be two signed compositions of $n$
and let $d \in X_{CD}$. Then $d^{-1} \in X_{DC}$ and, by
Proposition \ref{conjugaison x}, we have that
$$X_{C \cap \lexp{d}{D}} d = X_{\lexp{d^{-1}}{C} \cap D}.$$

\medskip

\begin{cor}\label{bijection double}
The map $X_{CD} \to W_C \backslash W_n /W_D$ is bijective.
\end{cor}

\begin{proof}
The proposition \ref{double classe} (f) shows that the map is
injective. The surjectivity follows from the fact that, if $w \in
W_n$ is an element of minimal length in $W_C w W_D$, then $w \in
X_{CD}$.
\end{proof}

\medskip

\begin{cor}\label{un cas facile}
If $C$ is parabolic or if $D$ is semi-positive, then
$$X_D=\coprod_{d \in X_{CD}} X_{C \cap \lexp{d}{D}}^C d.$$
\end{cor}

\begin{proof}
It follows from Corollary \ref{bijection double} that
$$|X_D|=|W_n/W_D|=\sum_{d \in X_{CD}} |W_C d W_D/W_D|=
\sum_{d \in X_{CD}} |X_{C \cap \lexp{d}{D}}|,$$ the last equality
following from Proposition \ref{double classe} (d). So, it remains
to show that, if $d \in X_{CD}$ and if $x \in X_{C \cap
\lexp{d}{D}}^C$, then $xd \in X_D$.

Assume that we have found $s \in S_D$ such that $\ell(xds) <
\ell(xd)$. If $s \in T_D$, then $s=t_i$ for som $i \in I_n$. But,
by Lemma \ref{easy}, $xd(i) > 0$, so $\ell(xdt_i) > \ell(xd)$,
contradicting our hypothesis. Therefore, $s \in S_{D^-}$, that is
$s=s_i$ for some $i \in [1,n-1]$. If $C$ is parabolic, $C \cap
\lexp{d}{D}$ is also parabolic. Therefore, $\ell(xds)>\ell(xd)$
which is a contradiction, so $D$ is semi-positive. Therefore, we
have that $t_i$
and $t_{i+1}$ belong to $T_D$. Thus, by Lemma \ref{easy}, we have
$xd(i) > 0$ and $xd(i+1) > 0$. Moreover, since $\ell(xds_i) <
\ell(xd)$, we have
$$0 < xd(i+1) < xd(i).\leqno{(*)}$$
But, since $d \in X_D$, we have $d(i+1) > d(i)$. So, by
Proposition \ref{double classe} (b), we have that $d s_i d^{-1}
\in S_{C \cap \lexp{d}{D}}$. Thus $\ell(x(ds_id^{-1})) > \ell(x)$
because $x \in X_{C \cap \lexp{d}{D}}^C$. In other words, $xd(i+1)
> xd(i)$. This contradicts $(*)$.
\end{proof}

If $E$ is a signed composition of $n$ such that $C \subset E$ and
$D \subset E$, we set $X_{CD}^E=X_{CD} \cap W_E$.

\medskip

\ex It is not true in general that $X_D=\coprod_{d \in X_{CD}}
X_{C \cap \lexp{d}{D}}^C d$. This is false, if
$n=k+l$ with $k$, $l \ge 1$, $C=(\bbar k,\bbar l)$ and $D=(n)$.
See Example~\ref{X union double} for precisions.

 In \cite{bbht}, the authors has given a proof of the Solomon theorem using
tools which sound like the above results. Here, we cannot  translate their proof because
 of the complexity  of the decomposition of $X_D$ (which involve negative coefficients).

\subsection{A partition of ${\boldsymbol{W_n}}$}
If $C=(c_1,\dots,c_r)$ is a signed composition of $n$, we set
$$A_C=\{s_{|c_1|+\dots+|c_i|}~|~i \in [1,r]\text{ and } c_i < 0
\text{ and } c_{i+1} > 0\}$$
$$\SC_C=S_C' \coprod A_C.\leqno{\text{and}}$$
As example, $A_{(1, \bar 3, \bar 1,2,\bar 1, 1)}=\{s_5 , s_8\}$.
Note that $\SC_C=\SC_D$ if and only if $C=D$. If $w \in W_n$, then
we define {\it the ascent set of $w$}:
$$\UC_n'(w)=\{s \in S_n'~|~\ell(ws) > \ell(w)\}.$$
Finally, following Mantaci-Reutenauer, we associate to each
element $w \in W_n$ a signed composition $\Cb(w)$ as follows.
First, let $\Cb^+(w)$ denote
the biggest composition (for the order $\subset$)
of $n$ such that, for every $1 \le i \le r$,
the map $w : I_{\Cb^+(w)}^{(i)} \to I_n$ is increasing and has
constant sign. Now, we define $\nu_i = \sgn(w(j))$ for
$j \in I_{\Cb^+(w)}^{(i)}$. The {\it descent composition} of $w$ is
$\bs C(w)=(\nu_1 c_1^+, \dots,\nu_r c_r^+)$.\\

\ex $\bs C (\underbrace{9. \bar 3 \bar 2\bar 1 . \bar 4. 58. \bar 6.
7}_{\in W_9}) = (1, \bar 3, \bar 1,2,\bar 1, 1) \VDash 9 $.\\

\noindent The following proposition is easy to check (see Remark
\ref{superieur}):

\begin{prop}\label{composition associee}
If $w \in W_n$, then $\UC_n'(w)=\SC_{\Cb(w)}$.
\end{prop}

\rem Mantaci and Reutenauer have defined the {\it descent shape} of a signed
permutation \cite{mantaci-reutenauer}. It is a signed composition defined
similarly than descent composition except that the absolute value of the
letters in $u_i$ must be in increasing order. For instance, the descent shape
of $9. \bar 3. \bar 2.\bar 1 \bar 4. 58. \bar 6.7$ is $(1,\bar 1, \bar 1, \bar 2 , 2 , \bar 1 , 1)$.\\

\medskip

\exemple{descente eta C}
Let $n'$ be a non-zero natural number, $n' < n$ and let $c \in \ZM$
such that $n-n'=|c|$. Let $w \in W_{n'} \subset W_n$ and
write $\Cb(w)=(c_1,\dots,c_r) \VDash n'$. Then
$\Cb(\eta_{(n',c)} w)=(c_1,\dots,c_r,c)$. Consequently, if $C \VDash n$,
an easy induction argument shows that $\Cb(\eta_C)=C$.

\medskip

We have then defined a surjective map
$$\Cb : W_n \longto \Comp(n)$$
whose fibers are equal to those of the application
$\UC_n' : W_n \to \PC(S_n')$. The surjectivity follows
from Example \ref{descente eta C}. If $C \VDash n$, we define
$$Y_C=\{w \in W_n~|~\Cb(w)=C\}.$$
Then
$$W_n=\coprod_{C \VDash n}Y_C.$$

\exemple{Y elementaire}
We have $Y_n=\{1_n\}$, $Y_\nba=\{\s_n w_n\}$,
$Y_{(1,\dots,1)}=\{\s_n\}$ and $Y_{(\bar 1 ,\dots ,\bar 1)}=\{w_n\}$.

\bigskip

First, note the following elementary facts.

\begin{lem}\label{pre ordre}
Let $C$ and $D$ be two signed compositions of $n$. Then:
\begin{itemize}
\itemth{a} If $Y_C \cap X_D \not= \varnothing$, then $Y_C \subset X_D$.

\itemth{b} $\eta_C \in Y_C$ and $Y_C \subset X_C$.
\end{itemize}
\end{lem}

\begin{proof}
(a) If $w \in W_n$, then $w \in X_D$ if and only if
$\UC_n'(w)$ contains $S_D'$. Since the map $w \mapsto \UC_n'(w)$
is constant on $Y_C$ (see Proposition \ref{composition associee}),
(a) follows.

\medskip

(b) By Example \ref{descente eta C}, we have $\eta_C \in Y_C \cap X_C$.
Therefore, by (a), $Y_C \subset X_C$.
\end{proof}

We then define a relation $\leftarrow$ between signed composition
of $n$ as follow. If $C$, $D \VDash n$, we write $C \leftarrow D$
if $Y_D \subset X_C$. We denote by $\infspe$ the transitive closure
of the relation $\leftarrow$.
It follows from Lemma \ref{pre ordre} (a) that
\equat\label{X union}
X_C=\coprod_{C \leftarrow D} Y_D.
\endequat

\exemple{x reunion y}
Let $w \in W_n$. By Remark \ref{superieur}, $w \in X_\nba$
if and only if the sequence $(w(1),w(2),\dots,w(n))$
of elements of $I_n$ is strictly increasing (see Remark~\ref{superieur}).
So there exists a unique $k \in \{0,1,2,\dots,n\}$ such that $w(i) > 0$ if
and only if $i > k$. Note that $k=\ell_t(w)$.
Let $i_1 < \dots < i_k$ be the sequence of elements
of $I_n$ such that $(w(1),\dots,w(k))=(\bar i_k,\dots, \bar i_1)$.
Then $w=r_{i_1} r_{i_2} \dots r_{i_k}$ where, if $1 \le i \le n$,
we set $r_i=s_{i-1} \dots s_2 s_1 t$. Note that $\Cb(w)=(\kba,n-k)$.
Therefore,
$$X_\nba=\{r_{i_1} r_{i_2} \dots r_{i_k}~|~
0 \le k \le n\text{ and } 1 \le i_1 < i_2 < \dots < i_k \le n\}.$$
Note that $\ell(r_{i_1} r_{i_2} \dots r_{i_k})=i_1+i_2 + \dots + i_k$
and $\ell_t(r_{i_1} r_{i_2} \dots r_{i_k})=k$. We get
$$X_\nba=\coprod_{0 \le k \le n} Y_{(\kba,n-k)},$$
and, for every $k \in \{0,1,2,\dots,n\}$, we have
$$Y_{(\kba,n-k)}=\{r_{i_1} r_{i_2} \dots r_{i_k}~|~
1 \le i_1 < i_2 < \dots < i_k \le n\}.$$
This shows that $(\nba) \leftarrow (\kba,n-k)$.

\bigskip

\begin{prop}\label{ordre}
Let $C$ and $D$ be two signed compositions of $n$. Then:
\begin{itemize}
\itemth{a} $C \leftarrow C$.

\itemth{b} If $C \subset D$, then $C \leftarrow D$.

\itemth{c} $\infspe$ is an order on $\Comp(n)$.
\end{itemize}
\end{prop}

\begin{proof}
(a) follows immediately from Lemma \ref{pre ordre} (b).

\medskip

(b) If $C \subset D$, then $X_D \subset X_C$. But, by Lemma
\ref{pre ordre} (b), we have $Y_D \subset X_D$.
So $C \leftarrow D$.

\medskip

(c) Let $a_C=\ell(\mu_C)$. By (a), $\infspe$ is reflexive. By
definition, it is transitive. So it is sufficient to show that it
is antisymmetric. But it follows from Lemma \ref{pre ordre} (b)
that:

\begin{quotation}
\noindent\quad $\bullet$ If $C \leftarrow D$, then $a_D \le a_C$.

\noindent\quad $\bullet$ If $C \leftarrow D$ and if $a_C=a_D$, then $C=D$.
\end{quotation}

\noindent The assertion (c) now follows easily from these two remarks.
\end{proof}

\bigskip

\exemple{X union double} If $C=(c_1,\dots,c_r)$ is  a composition of $n$ (not a signed composition), we will prove that
$$X_\nba=\coprod_{\begin{array}{c}
           \SS{0 \le m_2 \le c_2} \\
           \SS{0 \le m_3 \le c_3} \\
           {\dots} \\
           \SS{0 \le m_r \le c_r} \\
       \end{array}}
X_{(\cba_1,m_2,c_2-m_2,\dots,m_r,c_r-m_r)}^C
Y^{(\cba_1,m_2,c_2-m_2,\dots,m_r,c_r-m_r)}_{(\cba_1,
\mba_2,c_2-m_2,\dots,\mba_r,c_r-m_r)}
\s_{C,m_2,\dots,m_r}^{-1},$$
where $\s_{C,m_2,\dots,m_r} \in \SG_n$ satisfies
$$\s_{C,m_2,\dots,m_r}(S_{(c_1,m_2,c_2-m_2,\dots,m_r,c_r-m_r)}')
\subset S_n'$$
$$\s_{C,m_2,\dots,m_r} \in
X_{(c_1,m_2,c_2-m_2,\dots,m_r,c_r-m_r)}.\leqno{\text{and}}$$
%We want to prove that
%$$X_\nba=\coprod_{\begin{array}{c}
%           \SS{0 \le m_2 \le c_2} \\
%           \SS{0 \le m_3 \le c_3} \\
%           {\dots} \\
%           \SS{0 \le m_r \le c_r} \\
%      \end{array}}
%X_{(\cba_1,m_2,c_2-m_2,m_3,c_3-m_3,\dots,m_r,c_r-m_r)}^C
%Y^{(\cba_1,m_2,c_2-m_2,m_3,c_3-m_3,\dots,m_r,c_r-m_r)}_{(\cba_1,
%\mba_2,c_2-m_2,\mba_3,c_3-m_3,\dots,\mba_r,c_r-m_r)}
%\s_{C,m_2,\dots,m_r}^{-1},$$
%where $\s_{C,m_2,\dots,m_r} \in \SG_n$ satisfies
%$$\s_{C,m_2,\dots,m_r}(S_{(c_1,m_2,c_2-m_2,m_3,c_3-m_3,\dots,m_r,c_r-m_r)}')
%\subset S_n'$$
%$$\s_{C,m_2,\dots,m_r} \in
%X_{(c_1,m_2,c_2-m_2,m_3,c_3-m_3,\dots,m_r,c_r-m_r)}.\leqno{\text{and}}$$
By an easy induction argument, it is sufficient to prove it whenever
$r=2$. In other words, we want to prove that, if $k+l =n$ with
$k$, $l \ge 0$, then
$$X_\nba=\coprod_{0 \le m \le l} X_{(\kba,m,l-m)}^{(\kba,\lba)}
Y^{(\kba,m,l-m)}_{(\kba,\mba,l-m)} \s_{k,l,m}^{-1},\leqno{(*)}$$
where $\s_{k,l,m} \in \SG_n$ satisfies
$\s_{k,l,m}(S_{k,m,l-m}') \subset S_n'$ and $\s_{k,l,m} \in X_{(k,m,l-m)}$.
But, if $0 \le m \le l$, we set
$$\s_{k,l,m}(i)=\left\{\begin{array}{ll}
m+i & \text{ if } 1 \le i \le k, \\
i-k & \text{ if } k+1 \le i \le k+m, \\
i & \text{ if } k+m+1 \le i \le n, \\
\end{array}\right.$$
and one can easily check that $(*)$ holds. Moreover, since
$S_{k,m,l-m}'=S_n' \setminus \{s_k,s_{k+m}\}$, we get
that $\s_{k,l,m}(S_{k,m,l-m}') \subset S_n'$ and
$\s_{k,l,m} \in X_{(k,m,l-m)}$.

\medskip

%%%%%
%%
%% Generalized descent algebra\label{section solomon}
%%
%%%%%

\section{Generalized descent algebra\label{section solomon}}

\subsection{Definition}
If $C$ and $D$ are two signed compositions of $n$ such that
$C \subset D$, we set
$$x_C^D=\sum_{w \in X_C^D} w \qquad \in \ZM W_D$$
$$y_C^D=\sum_{w \in Y_C^D} w \qquad \in \ZM W_D.\leqno{\text{and}}$$
Now, let
$$\Sigma'(W_D)=\mathop{\oplus}_{C \subset D} \ZM y_C^D \qquad
\subset \ZM W_D.$$
Note that
$$\Sigma'(W_D)=\mathop{\oplus}_{C \subset D} \ZM x_C^D$$
by \ref{X union} and Proposition \ref{ordre}. We define
$$\th_D : \Sigma'(W_D) \longto \ZM \Irr W_D$$
as the unique $\ZM$-linear map such that
$$\th_D(x_C^D)=\Ind_{W_C}^{W_D} 1_C$$
for every $C \subset D$.
Here, $1_C$ is the trivial character of $W_C$. We denote by $\e_D$
the sign character of $W_D$.

\medskip

\noindent{\it Notation.}
If $D=(n)$, we set $x_C^D=x_C$, $y_C^D=y_C$ for
simplification. If $E$ is $\ZM$-module, we denote by $\QM E$
the $\QM$-vector space $\QM \otimes_\ZM E$. We denote by $\th_{D,\QM}$
the extension of $\th_D$ to $\QM\Sigma'(W_D)$ by $\QM$-linearity.

\medskip

\rem $\Sigma'(W_n)$ contains the Solomon descent algebras
of $W_n$ and $\Sn_n$. Moreover,  $\Sigma'(W_n)$ is precisely  the
{\it Mantaci-Reutenauer algebra} which is, by definition, generated by
$y_{D}=y_{D}^{(n)}$, for all $D \VDash n$.

\medskip

\subsection{First properties of ${\boldsymbol{\th_D}}$}
By the Mackey formula for product of induced characters and
by Proposition \ref{double classe}, we have
that
\equat\label{produit caracteres}
\th_n(x_C)\th_n(x_D)=\sum_{d \in X_{CD}}
\th_n(x_{\lexp{d^{-1}}{C} \cap D}).
\endequat

\medskip

\exemple{mackey facile}
If $C$ is parabolic or $D$ is semi-positive, then, by Corollary
\ref{un cas facile}, we have
$$x_D=\sum_{d \in X_{CD}} x_{C \cap \lexp{d}{D}}^C d.$$
Therefore, by Proposition \ref{simple classe} (b) and Remark \ref{xcd},
we get
$$x_C x_D = \sum_{d \in X_{CD}} x_{\lexp{d^{-1}}{C} \cap D}.$$
So $x_C x_D \in \Sigma'(W_n)$ and, by \ref{produit caracteres},
%\equat\label{theta parabolique}
$\th_n(x_C x_D)=\th_n(x_C)\th_n(x_D)$.\finl
%\endequat

\medskip

Before starting the proof of the fact that $\Sigma'(W_D)$
is a subalgebra of $\ZM W_D$ and that $\th_D$ is a morphism
of algebras, we need the following
result, which will be useful for arguing by induction.
If $C \subset D$, the transitivity of induction and Proposition
\ref{simple classe} (b) show that the diagram
\equat\label{theta induit}
\diagram
\Sigma'(W_C) \rrto^{\displaystyle{x_C^D .}}
\ddto_{\displaystyle{\th_C}} &&\Sigma'(W_D) \ddto_{\displaystyle{\th_D}} \\
&& \\
\ZM\Irr W_C \rrto^{\displaystyle{\Ind_{W_C}^{W_D}}} && \ZM \Irr W_D
\enddiagram
\endequat
is commutative.

Now, let $p_D : W_D \to \SG_{D^+}$ be the canonical projection.
It induces an injective morphism of $\ZM$-algebras
$p_D^* : \ZM\Irr \SG_{D^+} \to \ZM\Irr W_D$. Moreover,
the algebra $\Sigma'(\SG_{D^+})$ coincides with the usual
descent algebra in symmetric groups
and is contained in $\Sigma'(W_D)$.
Also, the diagram
\equat\label{diagramme sym}
\diagram
\Sigma'(\SG_{D^+}) \xto[0,2]|<\ahook \ddto_{\DS{\th_{D^-}}} &&
\Sigma'(W_D) \ddto^{\DS{\th_D}} \\
&& \\
\ZM\Irr \SG_{D^+} \xto[0,2]|<\ahook^{\DS{p_D^*}} && \ZM\Irr W_D \\
\enddiagram
\endequat
is commutative.

\medskip

\exemple{caracteres lineaires}
We have $y_{(\bar 1,\dots , \bar 1)}=w_n$, $y_\nba=w_n \s_n=\s_n w_n$, $y_{n}=1$
and $y_{(1,\dots,1)}=\s_n$. It is well-known \cite{solomon} that  $y_{(\bar 1,\dots , \bar 1)}$ belongs to
the classical descent algebra of $W_n$ and that
$$\th_n(w_n)=\e_n.\leqno{({\mathrm{a}})}$$
On the other hand,
$$\th_n(1_n)=1_{(n)}.\leqno{({\mathrm{b}})}$$
Also, by the commutativity of the diagram \ref{diagramme sym} and
as above, we have
$$\th_n(\s_n)=\g_n,\leqno{({\mathrm{c}})}$$
where $\g_n=p_n^* \e_\nba$. Finally, $w_n$ is a $\ZM$-linear
combination of $x_C$, where $C$ runs over the
parabolic compositions of $n$. Therefore, by Example \ref{mackey facile},
we have, for every $x \in \Sigma'(W_n)$,
$$\th_n(w_n x)=\th_n(w_n) \th_n(x) = \e_n \th_n(x).\leqno{({\mathrm{d}})}$$
In particular,
$$\th_n(y_\nba)=\th_n(w_n\s_n)=\e_n\g_n.\leqno{({\mathrm{e}})}$$
So we have obtained the four linear characters of $W_n$
as images by $\th_n$ of explicit elements of $\Sigma'(W_n)$.

\medskip

Let $\deg_D : \ZM\Irr W_D \to \ZM$ be the $\ZM$-linear
map sending an irreducible character of $W_D$ to its degree.
It is a morphism of $\ZM$-algebras.
Let $\aug_D : \ZM W_D \to \ZM$ be the augmentation morphism,
then it is clear that the diagram
\equat\label{augmentation}
\diagram
\Sigma'(W_D) \rrto^{\DS{\th_D}} \ddrrto_{\DS{\aug_D}} &&
\ZM\Irr W_D \ddto^{\DS{\deg_D}} \\
&&\\
&& \ZM
\enddiagram
\endequat
is commutative.

\bigskip

\subsection{Main result\label{main result}}
We are now ready to prove that $\Sigma'(W_D)$ is a $\ZM$-subalgebra
of $\ZM W_D$ and that $\th_D$ is a surjective morphism of algebras.

\bigskip

\begin{thm}\label{algebre}
Let $D$ be a signed composition of $n$,  then:
\begin{itemize}
\itemth{a} $\Sigma'(W_D)$ is a $\ZM$-subalgebra of $\ZM W_D$;

\itemth{b} $\th_D$ is a morphism of algebra;

\itemth{c} $\th_D$ is surjective and
$\Ker \th_D = \bigoplus_{\SS{C, C' \subset D} \atop \SS{C \equiv_D C'}}
\ZM(x_C^D - x_{C'}^D)$;

\itemth{d} $\Ker \th_{D,\QM}$ is the radical of the algebra
$\QM \Sigma'(W_D)$. Moreover,
$\QM \Sigma'(W_D)$ is a split algebra whose largest
semisimple quotient is commutative. In particular, all its simple
modules are of dimension $1$.
\end{itemize}
\end{thm}

\begin{proof}
We want to prove the theorem by induction on $|W_D|$. By taking
direct products, we may therefore assume that $D=(n)$ or
$D=(\nba)$. If $D=(\nba)$, then it is well-known that (a), (b), (c) and
(d) hold. So we may assume that $D=(n)$ and that
(a), (b), (c) and (d) hold for every signed composition $D'$
of $n$ different from $(n)$.

\medskip

(a) and (b):
Let $A$ and $B$ be two signed compositions of $n$. We want to prove
that $x_Ax_B \in \Sigma'(W_n)$ and that
$\th_n(x_A x_B)=\th_n(x_A)\th_n(x_B)$.
If $A$ is parabolic or $B$ is semi-positive, then
this is just Example \ref{mackey facile}. So we may assume that
$A$ is not parabolic and $B$ is not semi-positive.

First, note that $B \subset B^+$ and that $B^+$ is semi-positive.
Therefore, by Proposition~\ref{simple classe}(b) and Example \ref{mackey facile}, we have
$$x_A x_B= x_A x_{B^+} x^{B^+}_B = \sum_{d \in X_{A,B^+}}
x_{\lexp{d^{-1}}{A} \cap B^+} x^{B^+}_B = x_{B^+} \sum_{d \in X_{A,B^+}}
x_{\lexp{d^{-1}}{A} \cap B^+}^{B^+} x^{B^+}_B.$$
Assume first that $B^+\not=(n)$. Then, by induction hypothesis,
$$\sum_{d \in X_{A,B^+}}
x_{\lexp{d^{-1}}{A} \cap B^+}^{B^+} x^{B^+}_B \in \Sigma'(\Sn_{B^+})$$
and
$$\th_{B^+}\Bigl(\sum_{d \in X_{A,B^+}}
x_{\lexp{d^{-1}}{A} \cap B^+}^{B^+} x^{B^+}_B\Bigr)
=\sum_{d \in X_{A,B^+}} \th_{B^+}(x_{\lexp{d^{-1}}{A} \cap B^+}^{B^+})
\th_{B^+}(x^{B^+}_B).$$
Therefore, by \ref{theta induit} and \ref{diagramme sym}, $x_A x_B \in x_{B^+}\Sigma'(\Sn_{B^+})\subset\Sigma'(\Sn_{n})\subset  \Sigma'(W_n)$ and, by \ref{theta induit} and
by the Mackey formula for tensor product, we get
\begin{eqnarray*}
\th_n(x_A x_B) &=& \DS{
\Ind_{W_{B^+}}^{W_n} \Bigl(
\sum_{d \in X_{A,B^+}} \th_{B^+}(x_{\lexp{d^{-1}}{A} \cap B^+}^{B^+})
\th_{B^+}(x^{B^+}_B)\Bigr)} \\
&=& \DS{ \th_n(x_A)
\Ind_{W_{B^+}}^{W_n} \th_{B^+}(x^{B^+}_B) } \\
&=&\th_n(x_A)\th_n(x_B),
\end{eqnarray*}
as desired.

Therefore, it remains to consider the case where $B^+=(n)$. In particular,
$B=(n)$ or $(\nba)$. Since $B$ is not semi-positive, we have $B=(\nba)$.
By Example \ref{X union double}, we have
$$x_\nba=x_{A^-}^{A^+} +
\sum_{D \subset A^+} a_D x_D^{A^+} (\s_D^{-1}-1)$$
where $a_D \in \ZM$ and $\s_D(S_D') \subset S_n'$ and $\s_D \in X_D$
for every $D \subset A^+$. Therefore,
$$x_A x_B=x_A x_\nba = x_{A^+} \Bigl(x_A^{A^+} x_{A^-}^{A^+}
+ \sum_{D \subset A^+} a_D x_A^{A^+} x_D^{A^+} (\s_D^{-1}-1)
\Bigr).$$
Now, $W_A$ is a standard parabolic subgroup of $W_{A^+}$.
So, by Example \ref{mackey facile}, we have
$$x_A x_B=\sum_{d \in X_{A,A^-}^{A^+}}
x_{\lexp{d^{-1}}{A} \cap (A^-)}
+ \sum_{D \subset A^+} \Bigl(a_D \sum_{d \in X_{A,D}^{A^+}}
x_{\lexp{d^{-1}}{A} \cap D}(\s_D^{-1}-1)\Bigr).$$
Therefore, since $\s_D(S_D') \subset S_n'$ and $\s_D \in X_D$, we have that
$x_{\lexp{d^{-1}}{A} \cap D}\s_D^{-1}=
x_{\lexp{\s_D}{(\lexp{d^{-1}}{A \cap D})}}$. So $x_A x_B \in \Sigma'(W_n)$
and $\th_n(x_A x_B)=\th_n(x_A)\th_n(x_B)$ by the Mackey
formula for tensor product of induced characters. This concludes
the proof of (a) and (b). Indeed, the surjectivity of $\th_n$
is well-known.

\medskip

(c) First, let us show that $\th_n$ is surjective. Using the
induction hypothesis, the commutativity of the diagram \ref{theta induit},
and the classical description of irreducible characters of
$W_n$, we are reduced to prove that, for every $\chi \in \Irr \SG_n$,
$p_n^*(\chi)$ and $p_n^*(\chi)\e_n$ lie in the image of $\th_n$.
But it is well-known that $\th_\nba$ is surjective.
So the result follows from the commutativity of the diagram
\ref{diagramme sym} and from Example \ref{caracteres lineaires} (d).

Now, let $I=\sum_{\SS{C, C' \VDash n} \atop \SS{C \equiv_n C'}}
\ZM (x_C-x_{C'})$.
Then it is clear that $I \subset \Ker \th_n$.
Let $J=\oplus_{\l \in \Bip(n)} \ZM x_{\hat{\l}}$.
Then $\Sigma'(W_n)=I \oplus J$ and the map
$\th_n : J \to \ZM \Irr W_n$ is surjective. Since $J$ and
$\ZM \Irr W_n$ have the same rank (equal to $|\Bip(n)|$), we
get that $J \cap \Ker \th_n = 0$. So $I=\Ker \th_n$.

\medskip

(d) Let $R=\Rad(\QM\Sigma'(W_n))$ and
$K=\Ker \th_{n,\QM}$.
Since ${\text{Im}}(\th_{n,\QM})=\QM \Irr W_n$ is
a semisimple algebra, we get that $R \subset K$.

Now, let $\chi : \QM W_n \to \QM$ be the character of the $\QM W_n$-module
$\QM W_n$ (the regular representation). Then,
$\chi(w)=0$ for every $w \not= 1$.
Let $\chi'$ denote the restriction of $\chi$ to $\Sigma'(W_n)$.
We have $\chi'(x_C)=\chi(1)$ for every
$C \VDash n$. Therefore, $\chi'(x)=0$ for every $x \in K$
by (c). We fix now $x \in K$. Then, for every
$y \in \QM \Sigma'(W_n)$, we have $\chi'(xy)=0$ because
$xy \in K$ by (b). Since the
$\QM \Sigma'(W_n)$-module $\QM W_n$ is faithful, this
implies that $x \in R$. So $K \subset R$.
\end{proof}

\rem $\Sigma'(W_C) \simeq \Sigma'(W_{c_1}) \otimes \Sigma'(W_{c_2})\otimes\dots\otimes \Sigma'(W_{c_r})$.

\bigskip

\subsection{Further properties of ${\boldsymbol{\th_D}}$}
Let $\t_D : \ZM W_D \to \ZM$ be the unique linear map such that
$\t_D(w)=0$ if $w \not= 1$ and $\t_D(1)=1$. Then $\t_D$ is the canonical
symmetrizing form on $\ZM W_D$: in particular, the map
$\ZM W_D \times \ZM W_D \to \ZM$, $(x,y) \mapsto \t_D(xy)$ is a
non-degenerate symmetric bilinear form on $\ZM W_D$.
We denote by $\scalar{.}{.}_D$ the scalar product on
$\ZM\Irr W_D$ such that $\Irr W_D$ is an orthonormal basis.
The following property is a kind of ``isometry property'' for
the morphism $\th_D$.

\begin{prop}\label{isometrie}
If $x$, $y \in \Sigma'(W_D)$, then $\t_D(xy)=\scalar{\th_D(x)}{\th_D(y)}_D$.
\end{prop}

\begin{proof}
Let $C$ and $C'$ be two signed compositions of $n$ such that
$C$, $C' \subset D$. Then $\t_D(x_C^D x_{C'}^D)=|X_{CC'}^D|$ by definition of
$\t_D$. Moreover, since $\th_D(x_C)$ and $\th_D(x_{C'})$ take only rational
values, we have
$$\scalar{\th_D(x_C^D)}{\th_D(x_{C'}^D)}_D =
\scalar{\th_D(x_C^D)\th_D(x_{C'}^D)}{1_{W_D}}_D.$$
But, by \ref{produit caracteres} and by Frobenius reciprocity, we have
$$\scalar{\th_D(x_C^D)\th_D(x_{C'}^D)}{1_{W_D}}_D=|X_{CC'}^D|.$$
So the proposition follows now from the fact that $(x_C^D)_{C \subset D}$
generates $\Sigma'(W_D)$.
\end{proof}

\begin{cor}\label{ortho sigma}
$\Ker \th_D = \{x \in \Sigma'(W_D)~|~\forall~y \in \Sigma'(W_D),~\t_D(xy)=0\}$.
\end{cor}

\begin{proof}
Since $\scalar{.}{.}_D$ is non-degenerate on $\ZM\Irr W_D$,
this follows from Proposition \ref{isometrie}.
\end{proof}

\bigskip

Write $D=(d_1,\dots,d_r)$ and let 
$\Bip(D)$ denote the set of $r$-uples $(\l_{(1)},\dots,\l_{(r)})$
of bipartitions $\l_{(i)}=(\l_{(i)}^+,\l_{(i)}^-)$
such that $\l_{(i)}^-=\emptyset$ if $d_i < 0$
and $|\l_{(i)}|=|d_i|$ for every $i \in [1,r]$.
If $\l \in \Bip(D)$, we denote by $\CC_\l^D$ the conjugacy class
in $W_D$ of a Coxeter element of $W_{\hat{\l}}$ (with respect to
$S_{\hat{\l}}$). Let $f_\l^D$ denote the characteristic function of
$\CC_\l^D$. Then $f_\l$ is a primitive idempotent of $\QM\Irr W_n$.
Moreover, $(f_\l^D)_{\l \in \Bip(D)}$ is a complete family of orthogonal
primitive idempotents of
$\QM\Irr W_D$. Since $\th_D$ is surjective, there exists
a family of idempotents $(E_\l^D)_{\l \in \Bip(D)}$ of $\QM\Sigma'(W_D)$
such that

\begin{quotation}
\noindent (1) $\forall~\l \in \Bip(D)$, $\th_D(E_\l^D)=f_\l^D$.

\noindent (2) $\forall~\l, \m \in \Bip(D)$, $\l \not= \m \Rightarrow
E_\l^D E_\m^D=E_\m^D E_\l^D = 0$.

\noindent (3) $\DS{\sum_{\l \in \Bip(D)}} E_\l^D = 1$.
\end{quotation}

\bigskip

\begin{prop}\label{formule theta}
If $x \in \Sigma'(W_D)$, then
$$\th_D(x)=|W_D| \sum_{\l \in \Vdash D}
\frac{\t_D(xE_\l^D)}{|\CC_\l^D|} f_\l^D \in \ZM \Irr W_D.$$
\end{prop}

\begin{proof}
If $f \in \QM\Irr W_D$, then
$$f=|W_D| \sum_{\l \in \Bip(D)}
\frac{\scalar{f}{f_\l^D}_D}{|\CC_\l^D|} f_\l^D\in \ZM \Irr W_D.$$
If $f = \th_D(x)$ with $x \in \Sigma'(W_D)$, then we get the desired formula
just by applying Proposition \ref{isometrie} and the property (1) above.
\end{proof}

\medskip

\subsection{Character table\label{irreductible}}
Since all irreducible characters of $W_D$ have rational values,
the algebra $\QM\Irr W_D$ may be identified with the $\QM$-algebra
of central functions $W_D \to \QM$.
If $\l \in \Bip(D)$, we denote by $\eval_\l^D : \QM\Irr W_D \to \QM$,
$\chi \mapsto \chi(c_\l^D)$, where $c_\l^D$ is some element of $\CC_\l^D$
(for instance, a Coxeter element of $W_\lamh$).
Then $\eval_\l^D$ is a morphism of algebras: it is an irreducible
representation of $\QM\Irr W_D$. Moreover, $\{\eval_\l^D~|~\l \in \Bip(D)\}$
is a complete set of representatives of isomorphy classes of
irreducible representations of $\QM\Irr W_D$. Now, let $\QM_\l^D$
denote the $\QM\Sigma'(W_D)$-module whose underlying vector space
is $\QM$ and on which an element $x \in \QM\Sigma'(W_D)$ acts by
multiplication by $\pi_\l^D(x)=(\eval_\l^D \circ \th_D)(x)$.
Then, by Theorem \ref{algebre}, we get:

\begin{prop}\label{irr prop}
$\{\QM_\l^D~|~\l \in \Bip(D)\}$ is a complete set of
isomorphy classes of $\QM\Sigma'(W_D)$-modules. We have
$$\Irr(\QM\Sigma'(W_D))=\{\pi_\l^D~|~\l \in \Bip(D)\}.$$
\end{prop}

The {\it character table} of $\QM\Sigma'(W_D)$ is the square matrix
whose rows and the columns are indexed by $\Bip(D)$ and whose
$(\l,\m)$-entry is the value of the irreducible character
$\pi_\l^D(x_{\hat{\m}}^D)$. Note that
$$\pi_\l^D(x_{\hat{\mu}}^D)=\Bigl(\Ind_{W_{\hat{\mu}}}^{W_D} 1_{\hat{\mu}}\Bigr)
(c_\l^D).$$

\medskip

\noindent{\it Notation.}
If $D=(n)$, we denote $\CC_\l^D$, $f_\l^D$,
$E_\l^D$, $c_\l^D$, $\eval_\l^D$, $\QM_\l^D$ and $\pi_\l^D$
by $\CC_\l$, $f_\l$, $E_\l$, $c_\l$, $\eval_\l$, $\QM_\l$ and
$\pi_\l$ respectively.

\medskip

Now, if $\l$, $\m \in \Bip(D)$, we write $\l \subset \m$ if there
exists some $w \in W_D$ such that
$W_{\hat{\l}} \subset \lexp{w}{W_{\hat{\mu}}}$. By Proposition
\ref{conjugaison}, $\subset$ is a partial order on $\Bip(D)$.
For this partial order, the character table of $\QM\Sigma'(W_D)$ is
triangular~:

\medskip

\begin{prop}\label{triangularite}
If $\pi_\l^D(x_{\hat{\m}}^D) \not= 0$, then $\l \subset \m$.
\end{prop}

\medskip

\begin{proof}
We may, and we will, assume that $D=(n)$.
If $\pi_\l(x_{\hat{\mu}}) \not= 0$, then there exists $w \in W_n$ such
that $w c_\l w^{-1} \in W_{\hat{\mu}}$. Therefore, there exists
$\nu \in \Bip(\hat{\mu})$ and $w' \in W_{\hat{\mu}}$
such that $w'w c_\l w^{-1}w^{\prime-1}$ is a Coxeter element of
$W_{\hat{\nu}}$. Let $C$ denote the unique signed composition of $n$
such that $W_{\hat{\nu}}=W_C$ and let $\l'=\lamb(C)$. Then
$w'w c_\l w^{-1}w^{\prime-1}$ is conjugate to $c_{\l'}$. Therefore,
$\l=\l'$. This completes the proof of the proposition.
\end{proof}

\medskip

In the last section of this paper, we will give the character
table of $\Sigma'(W_2)$.
%
%
% Descent sets and Mantaci-Reutenauer algebra
%
%

%
% An algorithm
%

\subsection{Combinatorial description}\label{section refinement}

In $\Sn_n$, the refinement of compositions is useful to construct
$X_{C}$ from $Y_{D}$ without considering subsets of $S'_{\bar n}$.
The aim of this part is to describe such a procedure in our case.
 Start with an example,  consider $C = (\bar 2, 1)$, then the subsets of $S'_3$ containing $S_C= \set{s_1, t_3}$ are $\set{s_1 , s_2, t_3}=\SC_{(\bar 2, 1)}$; $\set{s_1, t_2,t_3 }=\SC_{(\bar 1,1 ,1)}$;
 $\set{s_1, t_1,t_2,t_3}=\SC_{(2, 1)}$, $\set{s_1,s_2,t_2,t_3}=\SC_{(\bar 1, 2)}$ and $S'_3 = \SC_{(3)}$.
 Observe  that $(1,2)$ (which corresponds to $\set{s_2,t_1,t_2,t_3}\not\supset S_C$) is not obtained. Here, we define a procedure which give $(\bar 2, 1)$, $(\bar 1, 1,
1)$, $(\bar 1, 2)$ and $(3)$, without to obtain $(1,2)$ .

\medskip

Let $C = (c_1 ,\dots ,c_k)\VDash n$, we write:
\begin{itemize}
\item $C\getsB D$ if $D=(a_1 ,b_1,a_2,b_2,\dots ,a_k,b_k)\VDash n$ such that for all $i\in [1,k]$ we have
$|a_i|+|b_i|=|c_i|$; $a_i=c_i$ (hence $b_i=0$) if $c_i>0$;  $a_i\leq 0\leq b_i$ if $c_i<0$  (remove the $0$ from the list $(a_1 ,b_1,a_2,b_2,\dots ,a_k,b_k)$).
That is, $D$ is obtained from $C$ by {\it broken negative parts operations}.
\item $C\getsR D$ if $C$ is finer that $D\VDash n$, that is, $D$ can be obtained from $C$ by summing consecutive parts of $C$ having the same sign ({\it refinement operations}).
\end{itemize}

\medskip

\exemple{ex refinement 1}  Let $C=(1, \bar 2,\bar 1)$, then
$$
\set{D\VDash 4 \st C\getsB D}=\set{(1, \bar 2,\bar 1),(1,\bar 1 ,1, \bar 1),(1 ,\bar 1, 1, 1),(1,2,\bar 1),(1, \bar 2, 1),(1,2,1)}.
$$

\medskip

\remarque{rem refinement} Let $C,D \VDash n$, then   we have
$C\gets D$ if and only if $S_C \subset \SC_D$.  We deduce easily
from definitions,  Lemma~\ref{pre ordre} and Example~\ref{x
reunion y} the following properties for any $i\in [1,k-1]$:
\begin{itemize}
\item if $\sgn c_i = \sgn c_{i+1}$, $C=(c_1,\dots ,c_i
,c_{i+1},\dots, c_k)\getsR (c_1,\dots ,c_i +c_{i+1},\dots, c_k)=D$,
and this  means that  $\SC_{D}=\SC_{C}\uplus\set{s_{|c_1|+\dots
+|c_i|}}$;

\item  if $c_i,c_{i+1}<0$, then
$$C=(c_1,\dots ,c_i , c_{i+1},\dots, c_k)\getsB
 (c_1,\dots , c_i +1, 1 ,c_{i+1},\dots, c_k)=D$$ (remove
the $0$ from the list), and this means that
$\SC_{D}=\SC_{C}\uplus\set{t_{|c_1|+\dots +|c_i|}}$;

\item if $c_i<0$ and $c_{i+1}>0$, then
 $$C=(c_1,\dots ,c_i , c_{i+1},\dots, c_k)\getsB
 (c_1,\dots , c_i +1, 1 ,c_{i+1},\dots, c_k)=D$$ (remove
the $0$ from the list), and this means that
$\SC_{D}\uplus\set{s_{|c_1|+\dots
+|c_i|}}=\SC_{C}\uplus\set{t_{|c_1|+\dots +|c_i|}}$. Moreover, as
$s_{|c_1|+\dots +|c_i|}\notin S_C$, we have $S_C \subset \SC_D$,
that is, $C\gets D$;

\item finally, if $c_i<0$, then 
$C\getsB (c_1,\dots ,c_i+1,1 ,c_{i+1},\dots, c_k)=D$ and  $$C
\getsB (c_1,\dots , c_i +2, 2 ,c_{i+1},\dots, c_k)=D'$$ (remove the
$0$ from the list), and this means that
$\SC_{D'}=\SC_{D}\uplus\set{t_{|c_1|+\dots +|c_i|}}$. Hence
$S_C\subset \SC_D \subset \SC_{D'}$.

\end{itemize}
In all these cases, we have $C\gets D$.

\medskip

\begin{thm}\label{refinement signed compositions} Let $C,D\Vdash n$, then $C\gets D$ if and only if there is $E\VDash n$ such that $C\getsB E \getsR D$. Moreover, $E$ is uniquely determined.
\end{thm}

\begin{proof} Suppose that  $E$ exists, then it is easy to check (using Remark~\ref{rem refinement} and induction) that $S_C \subset \SC_E\subset \SC_D$, which implies  $C\gets D$.

Now, suppose that $C\gets D$. As $S_C \subset \SC_D$,  it is easy to construct a unique $E\VDash n$ such that $\SC_E\cap T_n = \SC_D\cap T_n$ and $C \getsB E$ (hence $C \gets B$). It remains to show that $E\getsR D$, that is, to show that $\SC_E \cap S_\nba \subset  \SC_D\cap S_\nba$. Let $s_j \in \SC_E \cap S_\nba$, then  either $j\in [|c_1|+\dots +|c_{i-1}|+1,\dots , |c_1|+\dots +|c_{i}|-1]$ hence $s_j \in S_C \subset \SC_D$; or $j=|c_1|+\dots +|c_{i}|$ and $c_i <0$ and $c_{i+1}>0$ by definition. But refinement operations do not act on parts having not  the same sign, that is, $s_j\in \SC_D$.
\end{proof}

\exemple{}  Consider the signed composition $C = (1, \bar 2,\bar 1)$ . Then
 we obtain from Theorem~\ref{refinement signed compositions} and Example~\ref{ex refinement 1}
\begin{eqnarray*}
X_{(1, \bar 2,\bar 1)} &=& Y_{(1 ,\bar 2,\bar 1)} \cup Y_{(1 ,\bar 3)}
\cup\, Y_{(1,\bar 1 ,1, \bar 1)}
\cup\, Y_{(1 ,\bar 1, 1, 1)}\cup Y_{(1,\bar 1, 2)}
\cup\,  Y_{(1,2,\bar 1)}\cup Y_{(3,\bar 1)}  \\
&&\cup\, Y_{(1, \bar 2, 1)}
\cup\, Y_{(1,2,1)}\cup Y_{(3,1)}\cup Y_{(1,3)}\cup Y_{(4)} .
\end{eqnarray*}

\medskip

%%%%%
%%%%%
%
% Irreducible characters of $W_n$
%
%
%%%%%
%%%%%

\section{Coplactic space}\label{irreducible characters}

%As $\theta$ in Theorem \ref{algebre} is an epimorphism, a natural question is:
%what are the antecedents of the set of irreducible characters $\Irr (W_n)$ of $W_n$?
%The aim of this section is to construct such a set.
%
%In symmetric groups, J\H ollenbeck \cite{jollenbeck} has used the classical epimorphism
%of Solomon to construct $\Irr(\Sn_n)$. The fact that this construction can be
%generalized to $W_n$ suggests us that $\gsa$, with the epimorphism $\theta$,
%is the good tool to generalized other contructions (where the Solomon algebra appeared)
%to $W_n$ (see \cite{duchamp, loday-ronco, malvenuto-reutenauer}).
%
%First, we recall the {\it generalized Robinson-Schensted correspondence} in type $B$,
%with associated {\it Knuth} ({\it plactic}) and {\it dual-Knuth} ({\it coplactic})
%{\it relations}, introduced by Bonnaf\'e and Iancu in \cite{bonnafe-iancu}.
%
%In $\S$\ref{irreducible characters}, we use the generalized
%descent algebra, its epimorphism and the generalized
%Robinson-Schensted correspondence \cite{bonnafe-iancu} to give a
%new construction of the irreducible characters of $W_n$. Then, we define the
%{\it coplactic space} $\cop_n$, which can be viewed as a space of bitableaux
%which generalize the type $A$ tableaux Hopf algebra of Poirier-Reutenauer
%\cite{poirier-reutenauer}.
%
%Finally, we construct $\Irr (W_n)$ from $\cop_n$  for combinatorial reasons.
%The main tool using in the proof is an analog of the refinement of composition,
%which will be given in $\S$\ref{broken-refinement}.

\subsection{Robinson-Schensted correspondence for ${\boldsymbol{W_D}}$}
In \cite{stanley}, the author defined a bijection between $W_n$ and
a certain set of bitableaux, which sounds like a Robinson-Schensted
correspondence. Let us recall here some of his results.
A {\it bitableau} is a pair $T=(T^+,T^-)$ of tableaux. The {\it shape} of $T$
is the bipartition $(\l^+,\l^-)$, where $\l^+$ is the shape of $T^+$ and
$\l^-$ is the shape of $T^-$: it is denoted by $\sh T$. We note
$|T|=|\sh T|$. The bitableau
$T$ is said to be {\it standard} if the set of numbers in $T^+$ and $T^-$ is
$[1,m]$, where $m=|T|$, and if the fillings of $T^+$ and $T^-$ are increasing
in rows and in column.

Let $D \VDash n$. Write $D=(d_1,\dots,d_r)$ and denote by
$\SBT(D)$ the set of $r$-uples $T=(T_1,\dots,T_r)$ of
bitableaux $T_i=(T_i^+,T_i^-)$ such that $|T_i|=|d_i|$,
$T_i^-=\emptyset$ if $d_i < 0$, $T_i^+$ and $T_i^-$ are standard and
the fillings of $T_i^+$ and $T_i^-$ are exactly the numbers
in $I_{D,+}^{(i)}=[|d_1|+\dots+|d_{i-1}|+1,|d_1|+\dots+|d_i|]$.
The {\it shape} of $T$, denoted by $\sh T$,
is the $r$-uple of bipartitions $(\sh T_1,\dots,\sh T_r)$.
If $T \in \SBT(D)$, then $\sh T \in \Bip(D)$. If $\l \in \Bip(D)$,
we denote by $\SBT_\l^D$ the set of elements $T \in \SBT(D)$ such
that $\sh T=\l$. In \cite{stanley}, the author defined a bijection
(which we call
{\it generalized Robinson-Schensted correspondence})
% In \cite[Section 3]{bonnafe-iancu}, the authors have defined a bijection
\begin{eqnarray*}
\pi_D : W_D &\longto& \{(P,Q) \in \SBT(D) \times \SBT(D)~|~\sh P = \sh Q\} \\
w & \longmapsto & (\Pb_D(w),\Qb_D(w)).
\end{eqnarray*}
Note that, in \cite{stanley} (see also \cite[Section 3]{bonnafe-iancu}), the bijection has been defined
only for $D=(n)$. It is not difficult to deduce from this the bijection
$\pi_D$ for general $D$.
To this bijection is associated a partition of $W_n$ as follows:
if $Q \in \SBT(D)$, we set
$$Z_Q^D=\{w \in W_D~|~\Qb_D(w)=Q\}.$$
Then
$$W_D=\coprod_{Q \in \SBT(D)} Z_Q^D.$$

\bigskip

\subsection{Properties} 
First, note that the bijection $\pi_D$ satisfies the following property: 
if $w \in W_n$, then
\equat\label{inverse pi}
\pi_D(w^{-1})=(\Qb_D(w),\Pb_D(w)).
\endequat
In particular, if $Q$ and $Q'$ are two elements of $\SBT(D)$, then
\equat\label{inverse cardinal}
|Z_Q^D \cap (Z_{Q'}^D)^{-1}| = \begin{cases}
                               1 & \text{if } \sh Q =\sh Q',\\
                   0 & \text{otherwise.}
                   \end{cases}
\endequat

\medskip

\rem $\pi_\nba$ is the usual Robinson-Schensted correspondence.
For simplification, we denote by $Z_Q=Z_Q^{(n)}$ if $Q \in \SBT(n)$.

\medskip

 In \cite[Section 3]{bonnafe-iancu}, the authors give an another way to define the equivalence relation associated to this partition
which looks like coplactic equivalence or dual-Knuth equivalence.  If $w$, $w' \in W_D$, we
write $w \smile_D w'$ if $w'w^{-1} \in S_{D^-} \subset S_\nba=\{s_1,\dots,s_{n-1}\}$
and  $\DC_D'(w^{-1}) \not\subset \DC_D'(w^{\prime -1})$ and
$\DC_D'(w^{\prime -1}) \not\subset \DC_D'(w^{-1})$. Note that the relation
$\smile_D$ is symmetric. We denote by $\sim_D$ the reflexive and transitive
closure of $\smile_D$. It is an equivalence relation, called the
{\it coplactic equivalence relation}. The equivalence classes for this
relation are called the {\it coplactic classes} of $W_D$.
We denote by $\Cop(W_D)$ the set of coplactic classes for the relation
$\sim_D$.

By \cite[Proposition 3.8]{bonnafe-iancu}, we have, for every $w$, $w' \in W_D$,
\equat\label{knuth}
w \sim_D w' \Longleftrightarrow \Qb_D(w)=\Qb_D(w').
\endequat
So $\Cop(W_D)=\{Z_Q^D~|~Q \in \SBT(D)\}$.\\

\remarque{combinatorial coplactic} The relation $\smile_D$ has a useful
combinatorial interpretation (see
\cite[Proof of Proposition 3.8]{bonnafe-iancu}):
$w\smile_D s_iw$ ($s_i \in S_{D^-}$) if and only if:
\begin{itemize}
\item either $\sgn\big(w (i) \big)\not = \sgn\big(w (i+1) \big)$;
\item or $\sgn\big(w (i) \big) = \sgn\big(w (i+1) \big)$, then $s_i w$ is
obtained from $w$ by a classical dual-Knuth transformation; that is,
$w^{-1}(i-1)$ or $w^{-1}(i+2)$ lie between  $w^{-1}(i)$ and $w^{-1}(i+1)$.
\end{itemize}

\medskip

\begin{prop}\label{descente a droite}
Let $w,w' \in W_D$, then   $w \sim_D w' \implies \DC_D'(w)=\DC_D'(w')$.
\end{prop}

\begin{proof}
We may, and we will, assume that $D=(n)$.
If $T_0$ is a Young tableau, let
$\DC(T_0)=\{s_p \in S_{\bar{n}}~|~p+1$ lies in a row above
the row containing $p\}$. Let $\lexp{t}{T_0}$ denote the transposed tableau
of $T$. If $T=(T^+,T^-)$ is a standard bitableau, we set
$$\DC'(T) =\set{t_p  \st p\in T^-}\uplus
\set{s_p\st p\in T^+ \textrm{ and } p+1 \in T^-} \uplus \DC(T^+)
\uplus \DC\big(\lexp{t}{T^-}).$$ Then it is easy to check that
$\DC_n'(w)=\DC'(\Qb(w))$. This completes the proof of the
proposition.
\end{proof}

\medskip

\exemple{descent tableaux}
Let
$$
T =\Bigl(\quad \young(17,69,8)\quad,\quad \young(235,4)\quad\Bigr)
\quad \in \SBT(9).
$$
Then
\begin{eqnarray*}
\UD'(T)& =&\set{s_1 ,s_3, s_6  , s_8 , t_2 , t_3 ,t_4, t_5}.
\end{eqnarray*}

\medskip

\remarque{composition tableau} Using Proposition~\ref{descente a
droite}, we may assign a signed composition $\bs C(Q)\VDash n$ to
any standard bitableau $Q\in\SBT(n)$ by setting $\bs C(Q)=\bs
C(w)$ for any $w\in W_n$ such that $\bs Q(w)=Q$. One can determine
$\bs C(Q)$ directly from $Q$ thanks to the following procedure,
which is a combinatorial translation of the proof of
Proposition~\ref{descente a droite}. First one looks for maximal
subwords $j\;j+1\;j+2\ldots k$ of $1\;2\;3\ldots n$ such that
\begin{itemize}
\item either the numbers $j$, $j+1$, $j+2$, \dots, $k$ can be found in this
order in $Q^+$ when one goes from left to right (changes of rows are allowed)
\item or they can be found in this order in $Q^-$ when one goes from
top to bottom (changes of column are allowed).
\end{itemize}
The word $1\;2\;3\ldots n$ is then the concatenation of these
maximal subwords, and the signed composition $\bs C(Q)$ is the
sequence of the lengths of these subwords, adorned with a minus
sign if the letters of the subword can be found in $Q^-$. As an
example, consider $Q=(Q^+,Q^-)$ with
$$Q^+=\young(12678\thirteen,9\eleven\twelve,\ten)\qquad
\text{and}\qquad Q^-=\young(3\fourteen,4,5,\fifteen).$$ The
partition of $1\;2\dots 14\;15$ in maximal subwords is
$1\;2\mid3\;4\;5 \mid6\;7\;8\mid9\mid10\;11\;12\;13\mid14\;15$,
from what we can deduce that $\bs C(Q)=(2,\bar3,3,1,4,\bar2)$.

Therefore, we have by definitions
$$
X_C = \coprod_{C\gets \bs C(Q)} Z_Q . 
$$

\medskip

%%
%
%  Coplactic space
%
%%

\begin{prop}\label{coplacticite et induction}
Let $C,D \VDash n$ such that 
$C \subset D$. Let $w$, $w' \in W_C$ and $x$, $x' \in X_C^D$, then:
\begin{itemize}
\itemth{a} If $w \sim_C w'$, then $wx^{-1} \sim_D w'x^{-1}$.

\itemth{b} If $xw \sim_D x'w'$, then $w \sim_C w'$.

\itemth{c} If $w \sim_C w'$, then $w_C w \sim_C w_C w'$ and 
$w w_C \sim_C w' w_C$.
\end{itemize}
\end{prop}

\begin{proof}
(c) is clear. Let us now prove (a). We may assume that $w \smile_C w'$. 
But $\DC_C(w^{-1}) \subset \DC_D(x w^{-1})$ and 
$\DC_C(w^{\prime -1}) \subset \DC_D(x w^{\prime -1})$. So 
$wx^{-1} \smile_D w'x^{-1}$. 

\smallskip

We now prove (b). If $W_C$ is a standard parabolic subgroup of $W_D$, 
and using the fact that coplactic classes are left cells for a particular 
choice of parameters \cite[Theorem 7.7]{bonnafe-iancu}, 
then (b) follows from \cite{geck induction}. 
Therefore, by taking direct products and by arguing by induction on $|X_C^D|$, 
we may now assume that $D=(n)$ and $C=(k,l)$ with $k$, $l \ge 1$ and 
$k+l = n$. 

\smallskip

Let us start by proving (a). We may assume that $w \smile_C w'$. 
But $\DC_C(w^{-1}) \subset \DC_D(x w^{-1})$ and 
$\DC_C(w^{\prime -1}) \subset \DC_D(x w^{\prime -1})$. So 
$wx^{-1} \smile_D w'x^{-1}$. 

\medskip

Let us now prove (b). We may assume that $xw \smile_D x'w'$. 
Let $Q=\Qb_C(w)=\Qb_C(w')$. 
From Remark~\ref{combinatorial coplactic}, 
we have two cases: either $x'w'$ is obtained from $xw$ by a dual-Knuth 
relation, or $x'w'=s_i xw$ and $\sgn (xw(i))\not = \sgn (xw(i+1))$. 
In the first case, observe that, for any $i\in[1,k-1]$ and $i\in[k+1,k+l-1]$,   
$xw(k)<xw(k+1)$ if and only if $w(k)<w(k+1)$,  since 
$x\in X_{(k,l)}=X_{(\bar k,\bar l)}^{(\nba)}\subset W_\nba$. 
Then  we conclude by Remark~\ref{combinatorial coplactic} 
(which is exactly the result of Lascoux and Sch\"utzenberger on the 
shuffle of plactic classes \cite{lascoux-schutz}). 

In the second case, observe that, for any $k\in [1,n]$,
$$
(\star)\qquad \sgn (w(k))=\sgn (xw(k))=\sgn (s_i xw(k))
$$
since  $X_{(k,l)} = X_{(\bar k,\bar l)}^{(\nba)}\subset W_\nba$ and 
$s_i\in S_\nba$. If $s_i x = x'$, then $w=w'$ 
and the result follows. If $s_i x = x s_j$, with 
$s_j \in S_{(\bar k , \bar l)}$ (by Deodhar's Lemma), then 
$x'=x$.  Therefore $w'=s_j w$, by $(\star)$ and  
Remark~\ref{combinatorial coplactic}. So $w \smile_C w'$.
\end{proof}

\medskip

\subsection{Coplactic space}
Let $D \VDash n$. If $Q \in \SBT(D)$, we set
$$z_Q^D=\sum_{w \in Z_Q^D} w \qquad \in \ZM W_D.$$
Now, let
$$\QC_D=\mathop{\oplus}_{Q \in \SBT(D)} \ZM z_Q^D \qquad
\subset \ZM W_D$$
and
$$\QC_D^\perp=
\mathop{\oplus}_{\SS{Q, Q' \in \SBT(D)} \atop \SS{\sh Q=\sh Q'}}
\ZM(z_Q^D-z_{Q'}^D) \qquad
\subset \QC_D.$$
Then, by Proposition \ref{descente a droite}, we have
\equat
\Sigma'(W_D) \subset \QC_D.
\endequat
The next proposition justifies the notation $\QC_D^\perp$:

\begin{prop}\label{ortho}
$\QC_D^\perp=\{x \in \QC_D~|~\forall~y \in \QC_D,~\t_D(xy)=0\}$.
\end{prop}

\begin{proof}
Let $\QC_D'=\{x \in \QC_D~|~\forall~y \in \QC_D,~\t_D(xy)=0\}$.
Let $Q$ and $Q'$ be elements of $\SBT(D)$. Then, by \ref{inverse
cardinal}, we have \equat\label{scalaire z} \t_D(z_Q^D
z_{Q'}^D)=\begin{cases}
                1 & {\text{if }} \sh Q=\sh Q',\\
        0 & {\text{otherwise}}.
        \end{cases}
\endequat
This shows in particular that $\QC_D^\perp \subset \QC_D'$.

Let us now prove that $\QC_D' \subset \QC_D^\perp$.
Now, since $\QC_D/\QC_D^\perp$ is torsion free, it is sufficient
to prove that $\dim_\QM \QM\QC_D' \le \dim_\QM \QM\QC_D^\perp$.
But, by construction, we have
$\dim_\QM \QM\QC_D-\dim_\QM \QM\QC_D^\perp=|\Bip(D)|$.
Moreover, by Proposition \ref{isometrie}, we have
$\dim_\QM \QM\QC_D-\dim_\QM \QM\QC_D' \ge |\Irr W_D|=|\Bip(D)|$.
\end{proof}

\medskip

The next lemma is a generalization to our case of a result  of
Blessenohl and Schocker concerning the symmetric group
\cite{bless}.

\medskip

\begin{prop}\label{somme noyau}
We have $\QC_D=\Sigma'(W_D)+\QC_D^\perp$ and
$\Sigma'(W_D) \cap \QC_D^\perp = \Ker \th_D$.
\end{prop}

\medskip

\begin{proof}
Let us first prove that $\QC_D=\Sigma'(W_D)+\QC_D^\perp$. For this, we may,
and we will, assume that $D=(n)$. We first need to introduce an order on
bipartitions of $n$. We denote by $\le_{\slex}$ the lexicographic
order on $\Bip(n)$ induced by the following order on $I_n$:
$$\bar{1} <\slex \bar{2} <\slex \dots <\slex \bar{n} <\slex
1 <\slex 2 <\slex \dots <\slex n.$$
If $\l$ is a bipartition of $n$,
we denote by $Q_\l=Q(\eta_\lamh)$. If $\l=(\l^+,\l^-)$, then
it is easy to check that $\sh Q_\l=(\l^+,\lexp{t}{\l^-})=\l^*$, where
$\lexp{t}{\l^-}$ is the transpose of the partition $\l$, and, using Remarque~\ref{composition tableau}, that $Q_\l$ is obtained by numbered $Q_\l^+$ (resp. $\lexp{t}Q_\l^-$) the first column first, then the second one and so on. Now, let $Q \in \SBT(n)$. Then:

\medskip

\begin{quotation}
\begin{lem}\label{truc}
Assume that $Z_Q \subset X_\lamh$, then $\l \le_{\mathrm{sl}} (\sh Q)^*$.
Moreover, if $\sh Q=\l^*$, then $Q=Q_\l$.
\end{lem}

\begin{proof} First, we easily check (using Remark~\ref{composition tableau}), that  $\bs\l \big(\bs C(Q)\big)\le\slex (\sh Q)^*$ with equality if and only if $Q=Q_\l$.

Then, observe (using Theorem~\ref{refinement signed
compositions}),  that $\l\le\slex\bs\l\big(\bs C(Q)\big)$ with
equality if and only if $\bs C(Q)=\hat\l$. This conclude the
proof.
\end{proof}
\end{quotation}

\medskip

We are now ready to prove by descending induction on $(\sh Q)^*$ that
$z_Q \in \Sigma'(W_n)+\QC_n^\perp$.
If $(\sh Q)^*=(n,\emptyset)$, then
$Z_Q=\{1\}=Y_n$. So
$z_Q=y_n \in \Sigma'(W_n)$.

Now, assume that $(\sh Q)^* <\slex (n,\emptyset)$ and that $z_{Q'}
\in \Sigma'(W_n)+\QC_n^\perp$ for every $Q' \in \SBT(n)$ such that
$(\sh Q)^* <\slex (\sh Q')^*$. Let $\l=(\sh Q)^*$. Then
$z_Q=z_{Q_\l} + (z_Q-z_{Q_\l}) \in z_{Q_\l} + \QC_n^\perp$. So it
is sufficient to prove that $z_{Q_\l} \in \Sigma'(W_n) +
\QC_n^\perp$. But, by Lemma \ref{truc} $x_\lamh - z_{Q_\l}$ is a
sum of $z_{Q'}$ with $\l <_\lex (\sh Q)^*$. Hence, by the
induction hypothesis, we have $x_\l - z_{Q_\lamh} \in \Sigma'(W_n)
+ \QC_n^\perp$, as desired.

Now, let us prove that $\Sigma'(W_D) \cap \QC_D^\perp = \Ker \th_D$.
The natural map $\Sigma'(W_D) \to \QC_D/\QC_D^\perp$ is surjective,
so $\rank_\ZM \Sigma'(W_D) \cap \QC_D^\perp = \rank_\ZM \Ker \th_D$.
Since the $\ZM$-modules
$\Sigma'(W_D)/(\Sigma'(W_D) \cap \QC_D^\perp)$ and
$\Sigma'(W_D)/\Ker \th_D$ are torsion free, it is sufficient to
prove that $\Sigma'(W_D) \cap \QC_D^\perp$ is contained in
$\Ker \th_D$. But this follows from Proposition \ref{ortho} and
Corollary \ref{ortho sigma}.
\end{proof}

\bigskip

Using Proposition \ref{somme noyau}, we can easily extend the linear map
$\th_D$ to a linear map $\tilde{\th}_D : \QC_D \to \ZM\Irr W_D$.
If $x \in \QC_D$, write $x=a+b$ with $a \in \Sigma'(W_D)$ and
$b \in \QC_D^\perp$ and set
$$\tilde{\th}_D(x)=\th_D(a).$$
Then Proposition \ref{somme noyau} shows that $\tilde{\th}_D$ is
well-defined (that is, $\th_D(a)$ does not depend
on the choice of $a$ and $b$).

\medskip

\begin{thm}\label{theoreme theta}
Let $D\VDash n$. Then:
\begin{itemize}
\itemth{a} $\tilde{\th}_D$ is an extension of $\th_D$ to $\QC_D$;

\itemth{b} $\Ker \tilde{\th}_D=\QC_D^\perp$;

\itemth{c} if $x$ and $y$ are two elements of $\QC_D$, then
$\t_D(xy)=\scalar{\tilde{\th}_D(x)}{\tilde{\th}_D(y)}_D$;

\itemth{d} the diagram
$$\diagram
\QC_D \rrto^{\DS{\tilde{\th}_D}} \ddrrto_{\DS{\aug_D}} &&
\ZM\Irr W_D \ddto^{\DS{\deg_D}} \\
&&\\
&& \ZM
\enddiagram$$
is commutative;

\itemth{e} if $x \in \QC_D$, then
$$\tilde{\th}_D(x)= |W_D| \sum_{\l \in \Bip(D)}
\frac{\t_D(x E_\l^D)}{|\CC_\l^D|} f_\l^D.$$
\end{itemize}
\end{thm}

\begin{proof}
(a) and (b) are easy. (c) follows from Proposition \ref{isometrie}
and Proposition \ref{ortho}. (d) follows from the commutativity
of the diagram \ref{augmentation} and from the fact that
$\aug_D(\QC_D^\perp)=0$ (indeed, if $Q$ and $Q'$ are two elements
of $\SBT(D)$ of the same shape, then $|Z_Q^D|=|Z_{Q'}^D|$).
Using Proposition \ref{somme noyau}, it is sufficient to prove
(e) for $x \in \Sigma'(W_D)$ or $x \in \QC_D^\perp$.
If $x \in \Sigma'(W_D)$, this follows from
Proposition \ref{formule theta}. If $x \in \QC_D^\perp$,
this follows from Proposition \ref{ortho}.
\end{proof}

\medskip

\rem In the theorem, the case $D=(\bar n)$ is precisely the symmetric group case.

\medskip

\begin{cor}\label{unicite tilde}
If $\tilde{\th}$ is an extension of $\th_D$ to the $\QC_D$ such that, for all
$x$ and $y$ in $\QC_D$,
$\t_D(xy)=\scalar{\tilde{\th}(x)}{\tilde{\th}(y)}_D$, 
then $\tilde{\th}=\tilde{\th}_D$.
\end{cor}

\begin{proof}
Assume that $\tilde{\th}$ is an extension of $\th_D$ to $\QC_D$ such that
$\t_D(xy)=\scalar{\tilde{\th}(x)}{\tilde{\th}(y)}_D$ for all
$x$ and $y$ in $\QC_D$. Then, if $x \in \QC_D^\perp$ and
$\chi \in \ZM\Irr W_D$, then there exists $y \in \QC_D$ such that
$\tilde{\th}_D(y)=\chi$. So
$$\scalar{\chi}{\tilde{\th}(x)}_D=\scalar{\tilde{\th}(y)}{\tilde{\th}(x)}_D
= \t_D(xy)=0$$
by hypothesis and by Proposition \ref{ortho}. Since $\scalar{.}{.}_D$
is a perfect pairing on $\ZM \Irr W_D$, we get that
$\tilde{\th}(x)=0$.
So $\tilde{\th}$ coincides with $\tilde{\th}_D$ on $\Sigma'(W_D)$
and on $\QC_D^\perp$, so $\tilde{\th}=\tilde{\th}_D$
by Proposition \ref{somme noyau}.
\end{proof}

\medskip

Let $\l \in \Bip(D)$. Let $Q \in \SBT(D)$ be of shape $\l$.
Now, let
$$\x_\l=\tilde{\th}_D(z_Q).$$
Then $\x_\l$ depends only on $\l$ and not on the choice of $Q$.
Moreover, $\xi_\l \in \ZM\Irr W_D$, $\deg_D \x_\l = |Z_Q| > 0$
(see Theorem \ref{theoreme theta} (d))
and, by Theorem \ref{theoreme theta} (c) and \ref{scalaire z},
we have $\scalar{\xi_\l}{\xi_\l}_D=1$. This shows
that $\x_\l \in \Irr W_D$.
So we have proved the following proposition:

\begin{prop}\label{prop irred}
The map $\Bip(D) \to \Irr W_D$, $\l \mapsto \x_\l$ is
well-defined and bijective.
\end{prop}

\remarque{w0}
If $T=(T^+,T^-) \in \SBT(n)$, we denote by $T^\vee$ the standard
bitableau $(T^-,T^+)$. If $\l=(\l^+,\l^-) \in \Bip(n)$, we set
$\l^\vee=(\l^-,\l^+) \in \Bip(n)$. In particular, $\sh T^\vee=(\sh T)^\vee$.

Now, let $w \in W_n$.
Then $\pi_n(w_n w)=(\Pb(w)^\vee,\Qb(w)^\vee)$.
Therefore, if $Q \in \SBT(n)$ then $w_n Z_Q=Z_{Q^\vee}$. This shows
in particular that $w_n \QC_n = \QC_n$ and that
$w_n \QC_n^\perp=\QC_n^\perp$. Moreover,
\equat\label{w0 tilde}
\tilde{\th}_n(w_n z)=\e_n \tilde{\th}_n(z)
\endequat
for all $z \in \QC_n$. Indeed, this equality is true if
$z \in \Sigma'(W_n)$ by Theorem \ref{algebre} and it is obviously
true if $z \in \QC_n^\perp$. So we can conclude using Proposition
\ref{somme noyau}.
In particular, if $\l \in \Bip(n)$, then
\equat\label{w0 xi}
\xi_{\l^\vee}=\e_n \xi_\l.
\endequat

\medskip

\remarque{caractere symetrique}
Let $Q \in \SBT(n)$ be such that $Q^-=\emptyset$. Then $z_Q \in \QC_\nba$.
Therefore, $\QC_\nba \subset \QC_n$. Moreover,
$\QC_\nba^\perp = \QC_n^\perp \cap \QC_\nba$. Therefore, it follows
from the commutativity of Diagram \ref{diagramme sym} that the diagram
\equat\label{diagramme sym 2}
\diagram
\QC_\nba \xto[0,2]|<\ahook \ddto_{\DS{\tilde{\th}_\nba}} &&
\QC_n \ddto^{\DS{\tilde{\th}_n}} \\
&& \\
\ZM\Irr \SG_n \xto[0,2]|<\ahook^{\DS{p_n^*}} && \ZM\Irr W_n \\
\enddiagram
\endequat
is commutative. In particular, if $\l=(\l^+,\emptyset)$ is the
shape of $Q$, and if we denote by $\chi_{\l^+}^\nba$ the
irreducible character of $\SG_n$ associated to $\l^+$ (apply Proposition~\ref{prop irred} with $D=(\bar n)$), we have \equat\label{cas Sn} \xi_\l=p_n^*
\chi_{\l^+}^\nba.
\endequat

\bigskip

\subsection{Induction} 
We first start by an easy consequence of Proposition 
\ref{coplacticite et induction}. 

\medskip

\begin{lem}\label{induction coplaxique}
Let $C$, $D\VDash n$ be such  that $C \subset D$. Let $x \in \QC_C$. Then
\begin{itemize}
\itemth{a} $x_C^D x \in \QC_D$.

\itemth{b} If $x \in \QC_C^\perp$, then $x_C^D x \in \QC_D^\perp$.
\end{itemize}
\end{lem}

\begin{proof}
(a) By linearity, we may assume that $x=z_Q^C$ with $Q \in \SBT(C)$. 
Then, by Proposition \ref{coplacticite et induction} (a), 
we have that $X_C^D. Z_Q^C$ is a union of coplactic classes. So 
$x_C^D x \in \QC_D$. 

\medskip

(b) By linearity, we may assume that $x=z_Q^C-z_{Q'}^C$ where 
$Q$, $Q' \in \SBT(C)$ and $\sh(Q)=\sh(Q')$. We denote by 
$\psi : Z_Q^C \to Z_{Q'}^C$ the unique bijection such that 
$\Pb_C(\psi(w))=\Pb_C(w)$ for every $w \in Z_Q^C$. 

Then $x_C^D x = x_C^D .\sum_{w \in Z_Q^C} (w-\psi(w))$. 
But, if $a \in X_C^D$ and $w \in Z_Q^C$, then $\Pb_D(aw)=\Pb_D(a\psi(w))$ 
by Proposition \ref{coplacticite et induction} (b) and 
\ref{inverse pi}. We set $\psi'(aw)=a\psi(w)$, then the map 
$\psi' : X_C^D. Z_Q^C \to X_C^D.Z_{Q'}^C$ is bijective and satisfies 
$\sh \Qb_D(\psi'(w))=\sh \Qb_D(w)$ for every $w \in X_C^D. Z_Q^C$. 

Now, let $\l \in \Bip(D)$ and let $\EC_\l$ (resp. $\EC_\l'$) 
be the set of $w \in X_C^D. Z_Q^C$ (resp. $w \in X_C^D. Z_{Q'}^C$) 
such that $\sh \Qb_D(w)=\l$. Then 
$\psi'$ induces a bijection between $\EC_\l$ and $\EC_\l'$. 
Write $\EC_\l = \coprod_{i=1}^r Z_{Q_i}^D$ and 
$\EC_\l'=\coprod_{i=1}^{r'} Z_{Q_i'}^D$, using (a). Then, since 
$|Z_{Q_1}^D|=\dots=|Z_{Q_r}^D|=|Z_{Q_1'}^D|=\dots=|Z_{Q_{r'}}^D|$ and 
$|\EC_\l|=|\EC_\l'|$, we have $r=r'$. This shows that 
$x_C^D x \in \QC_D^\perp$.
\end{proof}

\medskip

\begin{cor}\label{induction Q}
Let $C$, $D\VDash n$ be such  that $C \subset D$. Then the 
diagram 
$$\diagram
\QC_C \rrto^{\displaystyle{x_C^D .}}
\ddto_{\displaystyle{\tilde{\th}_C}} &&
\QC_D \ddto_{\displaystyle{\tilde{\th}_D}} \\
&& \\
\ZM\Irr W_C \rrto^{\displaystyle{\Ind_{W_C}^{W_D}}} && \ZM \Irr W_D
\enddiagram$$
is commutative.
\end{cor}

\begin{proof}
This follows immediately from Proposition \ref{induction coplaxique} 
and from the commutativity of 
the diagram \ref{theta induit}.
\end{proof}

Now, if $\l \in \Bip(n)$, then we denote by $\chi_\l$ the irreducible 
character of $W_n$ associated to $\l$ via Clifford theory 
(see \cite{geck}). The link between the two parametrizations 
(the $\xi$'s and the $\chi$'s) is given by the following result:

\begin{cor}\label{explicite}
If $\l$ is a bipartition of $n$, then $\xi_\l=\chi_{\l^*}$.
\end{cor}

\begin{proof}
Write $\l=(\l^+,\l^-)$, $k=|\l^+|$ and $l=|\l^-|$. 
Let $Q^+$ be a standard tableau of shape $\l^+$ filled with 
$\{l+1,l+2,\dots,n\}$ and let $Q^-$ be a standard tableau of 
shape $\l^-$ filled with $\{1,2,\dots,l\}$. Then, 
by \cite[Proposition 4.8]{bonnafe-iancu}, 
$$Z_Q=X_{l,k} (w_l Z_{Q^-} \times Z_{Q^+}).$$
Therefore, by Corollary \ref{induction Q}, we have 
$$\xi_\l=\Ind_{W_{l,k}}^{W_n} (\tilde{\th_l}(w_l Z_{Q^-}^l) \boxtimes 
\tilde{\th_k}(Z_{Q^+}^k)).$$
So, by \ref{cas Sn} and by Remark \ref{w0}, we have 
$$\xi_\l=\Ind_{W_{k,l}}^{W_n} \Bigl(p_k^* \chi_{\l^+}^\kba \boxtimes 
\e_l (p_l^* \chi_{\l^-}^\lba)\Bigr).$$
The result now follows from \cite{geck}.
\end{proof}

%
%
% Hopf algebra of signed permutations
%
%

\section{Related Hopf algebras}

\medskip

\subsection{Hopf algebra of signed permutations}
Consider the graded $\ZM$-module
$$\SP={\mathop{\oplus}_{n\geq 0}} \ZM W_n,$$
where $W_0=1$. In \cite{aguiar-mahajan}, Aguiar and Mahajan have
shown that $\SP$ has a structure of Hopf algebra which is similar
to the structure of the Malvenuto-Reutenauer Hopf algebra on
permutations \cite{malvenuto-reutenauer}. Moreover, they have
shown that
$$
\Sigma' = \mathop{\oplus}_{n\geq 0} \Sigma'(W_n)
$$
is a Hopf subalgebra of $\SP$. We revise here the definition of
the product and the coproduct on $\SP$ with our point of view.

\medskip

\noindent{Notation - } If $C$ is a signed composition, then we
denote by $x_C$ the element of $\SP$ lying in $\ZM W_{|C|}$
corresponding to the $x_C$ defined in \S\ref{section solomon}.
Similarly, if $Q$ is a standard bitableau, then $z_Q \in \ZM
W_{|\sh Q|}$ is viewed as an element of $\SP$.

\medskip

Let $(u,v)\in W_n \times W_m$, we denote $u\times v$ the
corresponding element of $W_{n,m} \simeq W_n \times W_m$. If $w
\in W{n,m}$, we denote by $(w_{(n)}',w_{(m)}'')$ the corresponding
element of $W_n \times W_m$. We now define
$$u * v = x_{n,m}(u \times v) \in \ZM W_{n+m}.$$
We extend $*$ by linearity to a bilinear map $\SP \times \SP \to
\SP$.

Now, let $w \in W_n$. Then, for each $i \in [0,n]$, we denote by
$\pi_i(w)$ the unique element of $W_{i,n-i}$ such that $w \in
\pi_i(w) X_{i,n-i}^{-1}$. We set
$$\D(w) = \sum_{i=0}^n \pi_i(w)_{(i)}' \otimes_\ZM \pi_i(w)_{(n-i)}''
\in \SP \otimes \SP.$$ We extend $\D$ by linearity to a map $\D :
\SP \to \SP \otimes_\ZM \SP$.

\medskip

\remarque{combinatorial product} Combinatorially, we see this product as
follows: let $w=w_1\dots w_n$ be a  word of length $n$ in the alphabet
$I_n$, the {\it standardsigned permutation}  is the unique element
$\sts (w)  \in W_n$ such that
$$
\left\{\begin{array}{l}
\sts(w)(i) <\sts(w)(j)\iff\big(w_i<w_j) \quad\textrm{or}\quad
(w_i=w_j \textrm{ and } i<j  \big)\\
\textrm{and}\quad\sgn\big(\sts(w)(i)\big)=\sgn(w_i).
\end{array}\right.
$$
Then
$$
u\SPprod v = \sum_{w,w'} ww'
$$
where $ww'$ is the concatenation of $w$ and $w'$; and  the sum is
taken over all  words $w,w'$ on the alphabet $I_n$ such that
$\sts(w)=u$, $\sts(w')=v$ and  $\alph(u)\uplus\alph(v)=\entier{n}$
(where $\alph(u)=$ the set of absolute values of the letters in $u$).
For instance,  $\bar 1 2 \times 2 \bar 1 = \bar 1 2 4 \bar 3$ and
\begin{eqnarray*}
x_{(2,2)} &= &y_{(2,2)}+y_{(4)}\\
&=&1234+ 1324+1423+2314+2413+3412 .
\end{eqnarray*}
Hence
$
\bar 1 2 \SPprod 2 \bar 1 =\bar 1 2 4 \bar 3+ \bar1 3 4 \bar 2+
\bar 143\bar 2+\bar
2 3 4 \bar 1 + \bar 243\bar 1+\bar 342\bar 1.
$

\medskip

\remarque{combinatorial coproduct}
For $w\in W_n$ seen as a  word on the alphabet $I_n$ and
$i<j\in\entier n$, we denote  $w|\Dentier{i}{j}$ the subword
obtained by taking only the digits such that their absolute
values are in $\Dentier{i}{j}$.  Then we see combinatorially
the coproduct as
$$
\SPcoprod (w) = \sum_{i=0}^n w|\entier i \otimes
\sts\big(w|\Dentier{i+1}{n}\big) .
$$
As example,  consider $w=\bar2 3 1 \bar 4$, then we have the following
decompositions:
$$
w^{-1}= 3\bar1 2 \bar 4 = 3124(1\times \bar 12\bar 3)=1324 ( 2 \bar
1\times 1 \bar 2) = 1234(3\bar 12 \times \bar 1 ).
$$
Hence
$$
w=\bar2 3 1 \bar 4 = (1\times \bar 12\bar 3) 2314= (  \bar2
1\times 1 \bar 2)1324 = (\bar2 31 \times \bar 1 ) 1234.
$$
Thus
$$
 \SPcoprod (\bar 2 3 1 \bar 4)= \emptyset \otimes \bar 2 3 1
\bar 4+ 1 \otimes \bar 12  \bar 3 + \bar 2 1\otimes 1\bar 2 + \bar
23 1 \otimes \bar 1+ \bar 2 31 \bar 4 \otimes\emptyset .
$$

\medskip

\exemple{produit x} Let $C$ and $D$ be two signed composition. We
denote by $C \sqcup D$ the signed composition obtained by
concatenation of $C$ and $D$. Then
$$x_C * x_D = x_{C \sqcup D}.$$

\medskip

\exemple{coproduit x} We have
$$\D(x_n)=\sum_{i=0}^n x_i \otimes_\ZM x_{n-i}$$
$$\D(x_\nba)=\sum_{i=0}^n x_\iba
\otimes_\ZM x_{\overline{n-i}}.\leqno{\text{and}}$$

\medskip

We state here a result of Aguiar and Mahajan~\cite{aguiar-mahajan},
with our basis consisting of the $x_C$.

\begin{thm}\label{mesquin}
The graded vector space $\SP$, with the product $\SPprod$ and the
coproduct $\SPcoprod$ is a connected graded Hopf algebra; and
$\Sigma'$ is a Hopf subalgebra of $\SP$ which is freely generated
by elements $(x_n)_{n \in \ZM\setminus \{0\}}$ as algebra.
\end{thm}

\medskip

If $x$, $y \in \SP$, we define the product $xy \in \SP$ as follows:
if $x \in \ZM W_n$ and $y \in \ZM W_m$, then $xy=0$ if $m\not=
n$ and $xy$ coincides with the usual product $xy$ in $\ZM W_n$ if
$m=n$. Let $\t : \SP \to \ZM$ be the unique linear map which
coincides with $\t_n$ on $\ZM W_n$. The the map $\SP \times \SP \to
\ZM$, $(x,y) \mapsto \t(xy)$ is a scalar product on $\SP$. If $x$,
$y \in \SP$, we set
$$\t_\otimes(x \otimes y) = \t(x) \t(y).$$
The following proposition is easily checked from definitions:

\begin{prop}\label{self dual}
$\SP$ is self-dual for $\t$, that is,
$$
\t_\otimes\big((u\otimes v)\Delta(w)\big)=\t\big((u * v)w\big)
$$
for all $u$, $v$, $w \in \SP$.
\end{prop}

\medskip

\subsection{The  Hopf algebra of characters} We give here a
short recall of a result of Geissinger \cite{geissinger}. Consider
the graded $\ZM$-module
$$
\Hclass = \mathop{\oplus}_{n\geq 0} \ZM\Irr W_n .
$$
If $k$ and $l$ are two natural numbers, we denote by $\iota_{k,l}$
the canonical isomorphism
$$\iota_{k,l} : \ZM\Irr W_k \otimes_\ZM \ZM\Irr W_l \longmapright{\sim}
\ZM\Irr W_{k,l}.$$ Let $(\chi,\psi)\in\ZM\Irr W_k \times \ZM\Irr
W_l$. We define
$$
\chi \bullet \psi =
\Ind_{W_{k,l}}^{W_{k+l}}\iota_{k,l}\big(\chi\otimes_\ZM \psi\big)
\quad \in \ZM\Irr W_{k+l}.
$$
Now, let $\chi\in \class W_n$. We define
$$
\Res(\chi)= \sum_{i=0}^n  \iota_{i,n-i}^{-1}
\Res_{W_{i,n-i}}^{W_n} \chi \quad \in \mathop{\oplus}_{i=0}^n \ZM
\Irr W_i \otimes_\ZM \ZM\Irr W_{n-i} \subset \Hclass \otimes_\ZM
\Hclass.
$$
We denote $\scalar{\cdot}{\cdot}$ the unique scalar product on
$\Hclass$ which coincides with $\scalar{\cdot}{\cdot}_n$ on
$\ZM\Irr W_n$ and such that $\ZM\Irr W_n$ and $\ZM\Irr W_m$ are
orthogonal if $m \not= n$. We now define
$\scalar{\cdot}{\cdot}_\otimes$ on $\Hclass \otimes_\ZM \Hclass$
as follows: if $\chi$, $\chi'$, $\psi$, $\psi' \in \Hclass$, we set
$$
\scalar{\chi\otimes \psi}{\chi'\otimes
\psi'}_{\otimes}=\scalar{\chi}{\chi'} \scalar{\psi}{\psi'}.
$$

Geissinger \cite{geissinger} has shown that $\Hclass$ with product
$\bullet$ and
coproduct $\Res$ is a connected graded Hopf algebra. Moreover, for
any $\chi,\psi,\zeta\in\Hclass$,  the reciprocity law of Frobenius
can be viewed as \equat\label{propriete scalaire} \scalar{\chi
\otimes \psi}{\Res \zeta}_\otimes= \scalar{\chi \bullet
\psi}{\zeta}.
\endequat

\medskip

%
% The coplactic algebra
%

\subsection{The coplactic algebra and an Hopf epimorphism}

Let us intoduce
$$
\cop =\mathop{\oplus}_{n\geq 0} \cop_n.
$$
and
$$
\Hclass = \mathop{\oplus}_{n\geq 0} \ZM\Irr W_n .
$$
 We define $\th : \Sigma' \to \Hclass$ and
$\tilde{\th} : \QC \to \Hclass$ by
$$\th=\mathop{\oplus}_{n \ge 0} \th_n\qquad\text{and}\qquad
\tilde{\th}=\mathop{\oplus}_{n \ge 0} \tilde{\th}_n.$$

The first part of the following theorem shows that $\QC$ is 
a generalization of the Poirier-Reutenauer Hopf algebra of tableaux
\cite{poirier-reutenauer} to our case (see also \cite{bless}), and
the second part shows that J\"ollenbeck's 
construction generalizes to our case.

\begin{thm}\label{hopf_coplac}
$\QC$ is a Hopf subalgebra of $\SP$ containing $\Sigma'$.
 Moreover, $\th : \Sigma' \to \Hclass$ and $\tilde{\th} : \QC \to \Hclass$ are surjective Hopf algebra homomorphisms.
\end{thm}

\begin{proof} 
The fact that $\QC$ is a subalgebra of $\SP$ follows from
Proposition \ref{induction coplaxique}. 
To prove that it is a subcoalgebra, we proceed  as in the 
proof of  the result of Poirier and Reutenauer \cite{poirier-reutenauer}, 
using Remark~\ref{combinatorial coplactic}: let $Z$ be a coplactic class 
in $W_n$, $i\in [0,n]$ and $w\in Z$. Write $w= \pi_i (w) x$, where 
$x^{-1}\in X_{i,n-i}$. Let $u  \in W_i$ such that 
$u\smile_{(i)} \pi_i (w)'_{(i)}$. As 
$\sgn\big( x^{-1} w^{-1} (k)\big)=\sgn\big( w^{-1}(k)\big)$ and  
$x^{-1}(l)<x^{-1}(l+1)$, for all $l\in [1,i-1]$ and for all 
$l\in [i+1,n-1]$, we easily check that 
$(u\times \pi_i (w)''_{(n-i)})x \smile_{(n)} w$, using 
Remark~\ref{combinatorial coplactic}.  Let $v  \in W_{n-i}$ 
such that $v\smile_{(n-i)} \pi_i (w)''_{(n-i)}$, 
then $(u\times v )x \smile_{(n)} w$  as above. Therefore
$$
\SPcoprod \big(\sum_{w\in Z} w\big) = \sum_{i=0}^n
\sum_{Z_i , Z_{n-i}} \big(\sum_{u\in Z_i} u
\big) \otimes \big( \sum_{v\in Z_{n-i}} v \big),
$$
where $Z_i$ (resp. $Z_{n-i}$) are coplactic
classes in $W_i$ (resp. $W_{n-i}$).\\

We now need to prove that
$\tilde{\th}$ is an homomorphism of Hopf algebras. We first need a
lemma concerning the symmetric bilinear form $\b :
(\QC\otimes_\ZM\QC) \times (\QC \otimes_\ZM\QC) \to \ZM$, $(x,y)
\mapsto \t_\otimes(xy)$. Let 
$\tilde{\th}_\otimes =\tilde{\th} \otimes_\ZM  \tilde{\th} : 
\QC \otimes_\ZM \QC \to \Hclass \otimes_\ZM \Hclass$. Then:

\begin{quotation}
\begin{lem}\label{lemme noyau}
$\Ker \tilde{\th}_\otimes=\QC \otimes_\ZM \Ker \tilde{\th} + 
\Ker \tilde{\th} \otimes_\ZM \QC$ is the kernel of $\b$.
\end{lem}

\begin{proof}
By Theorem\ref{theoreme theta} (c), we have 
$$\b(x,y)=\scalar{\tilde{\th}_\otimes(x)}{\tilde{\th}_\otimes(y)}_\otimes$$
for all $x$, $y \in \QC \otimes_\ZM \QC$. This proves the lemma. 
\end{proof}
\end{quotation}

\medskip

Proposition \ref{self dual} and Lemma \ref{lemme noyau} show that
$\Ker \tilde{\th}$ is an ideal and a coideal of $\QC$. 
Since $\QC = \Sigma' + \Ker \tilde{\th}$, it is sufficient to 
prove that $\th$ is a bialgebra homomorphism. First, 
it is clear that $\tilde{\th}(x_C * x_D)=\tilde(x_{C \sqcup D})$. 
So $\tilde{\th}$ is an algebra homomorphism. Using this last 
propoerty and Theorem \ref{mesquin}, it is sufficient to prove 
that $\tilde{\th}_\otimes(\D(x_n))=\Res(\tilde{\th}(x_n))$ 
and $\tilde{\th}_\otimes(\D(x_{\bar{n}}))=\Res(\tilde{\th}(x_{\bar{n}}))$ 
But this follows easily from Example \ref{coproduit x}.
\end{proof}

\medskip

%%%%%
%%%%%
%%
%%  Example of $W_2$
%%
%%%%%
%%%%%

\section{The case $n=2$}

\medskip

In this Section, we will give a complete description of the
algebra $\Sigma'(W_2)$. For simplification, we set
$s=s_1$. Note that $t_1=t$
and $t_2=sts$. In other words, $S_2'=\{s,t,sts\}$.
Table I gives the correspondence between
reduced decomposition of elements of $W_2$ and permutations
of $I_2$ (if $w \in W_2$, we only give the couple $(w(1),w(2))$
since it determines $w$ as a permutation of $I_2$).
It also gives the value of $\UC_2'(w)$ and $\Cb(w)$.
Table II gives representatives of the conjugacy classes of
$W_2$. Table III gives, for each signed composition $C$ of $2$, the subgroup
$W_C$ of $W_2$, the set $S_C$,
the elements $x_C$ and $y_C$ of $\ZM W_2$ and also
gives the value of $\SC_C$.
Table IV provides
the decomposition of the induced characters
$\Ind_{W_{\lamh}}^{W_2} 1_{\lamh}=\th_2(x_\lamh)$ as a combination of the
$\xi_\m$, for $\l \Vdash 2$. Table V gives the character table
of $\QM\Sigma'(W_2)$ (see Subsection \ref{irreductible}).
We give in Table VI a complete set of orthogonal primitive
idempotents of $\QM\Sigma'(W_2)$.
Table VII gives the Cartan matrix of $\Sigma'(W_2)$. As usual,
the dots in the tables represent the number $0$.
Note that
$$w_2=stst=tsts.$$

We conclude the Section by a description of the algebra
$\QM\Sigma'(W_2)$ as a product of classical
indecomposable algebras.

\medskip

\noindent{\it Convention. } For avoiding the use of too many 
parenthesis, we have denoted by $\xi_\lamh$, $\pi_\lamh$ or $E_\lamh$ 
the objects $\xi_\l$, $\pi_\l$ or $E_\l$ respectively. For instance, 
$\xi_{1,\unbar}=\xi_{((1),(1))}$ and 
$\pi_{2}=\pi_{((2),\emptyset)}$ and $E_{\unbar,\unbar}=E_{(\emptyset,(1,1))}$.

\bigskip

\begin{centerline}
{\begin{tabular}{c}
\begin{tabular}{|c|c|c|c|}
\hline
\vertical $w$ & $(w(1),w(2))$ & $\UC_2'(w)$ & $\Cb(w)$ \\
\hline
\hline
\vertical $1$   & $(1      ,2)      $ & $\{s,t,sts\}$ & $(2)$               \\
\vertical $s$   & $(2      ,1)      $ & $\{t,sts\}$   & $(1,1)$             \\
\vertical $t$   & $(\bar{1},2)      $ & $\{s,sts\}$   & $(\bar{1},1)$       \\
\vertical $st$  & $(\bar{2},1)      $ & $\{s,sts\}$   & $(\bar{1},1)$       \\
\vertical $ts$  & $(2      ,\bar{1})$ & $\{t\}$       & $(1,\bar{1})$       \\
\vertical $sts$ & $(1      ,\bar{2})$ & $\{t\}$       & $(1,\bar{1})$       \\
\vertical $tst$ & $(\bar{2},\bar{1})$ & $\{s\}$       & $(\bar{2})$         \\
\vertical $w_2$ & $(\bar{1},\bar{2})$ & $\emptyset$   & $(\bar{1},\bar{1})$ \\
\hline
\end{tabular} \\
\\
\vertical {\bf Table I. Elements} \\
\end{tabular}
\qquad
\begin{tabular}{c}
\begin{tabular}{|c|c|c|}
\hline
\vertical $\lamh$           & $c_\l$ & $|\CC_\l|$ \\
\hline
\vertical $(2)$             & $st$   & $2$        \\
\vertical $(1,1)$           & $w_2$  & $1$        \\
\vertical $(1,\unbar)$      & $t$    & $2$        \\
\vertical $(\deuxbar)$      & $s$    & $2$        \\
\vertical $(\unbar,\unbar)$ & $1$    & $1$        \\
\hline
\end{tabular}\\
\\
\vertical {\bf Table II. Conjugacy classes} \\
\end{tabular}}
\end{centerline}

\bigskip

\begin{centerline}
{\begin{tabular}{c}
\begin{tabular}{|c|c|c|c|c|c|}
\hline
\vertical $C$ & $W_C$ & $S_C$ & $x_C$ & $y_C$ & $\SC_C$\\
\hline
\hline
\vertical $(2)$       & $W_2$      & $\{s,t\}$ & $1$          & $1$    & $\{s,t,sts\}$ \\
\vertical $(1,1)$     & $W_1 \times W_1$ & $\{t,sts\}$ & $1+s$ & $s$   & $\{t,sts\}$ \\
\vertical $(\bar{1},1)$ & $\SG_1 \times W_1$ & $\{sts\}$ & $1+s+t+st$ & $t+st$ & $\{t\}$ \\
\vertical $(1,\bar{1})$ & $W_1 \times \SG_1$ & $\{t\}$ & $1+s+ts+sts$ & $ts+sts$ & $\{s,sts\}$ \\
\vertical $(\bar{2})$ & $\SG_2$    & $\{s\}$   & $1+t+st+tst$ & $tst$   & $\{s,sts\}$ \\
\vertical $(\bar{1},\bar{1})$ & $1$ & $\emptyset$ & $\sum_{w \in W_2} w$ & $w_2$ & $\emptyset$ \\
\hline
\end{tabular}\\
\\
\vertical {\bf Table III. Bases of ${\boldsymbol{\Sigma'(W_2)}}$} \\
\end{tabular}}
\end{centerline}

\bigskip

\begin{centerline}
{\begin{tabular}{c}
\begin{tabular}{|c||c|c|c|c|c|}
\hline
\vertical & $\xi_2$ & $\xi_{1,1}$ & $\xi_{1,\unbar}$ & $\xi_\deuxbar$ & $\xi_{\unbar,\unbar}$ \\
\hline
\hline
\vertical  $\th_2(x_2)$               & 1 & . & . & . & . \\
\vertical  $\th_2(x_{1,1})$           & 1 & 1 & . & . & . \\
\vertical  $\th_2(x_{1,\unbar})$      & 1 & 1 & 1 & . & . \\
\vertical  $\th_2(x_\deuxbar)$        & 1 & . & 1 & 1 & . \\
\vertical  $\th_2(x_{\unbar,\unbar})$ & 1 & 1 & 2 & 1 & 1 \\
\hline
\end{tabular}\\
\\
\vertical {\bf Table IV. Decomposition of induced characters}\\
\end{tabular}}
\end{centerline}

\bigskip

\begin{centerline}
{\begin{tabular}{c}
\begin{tabular}{|c||c|c|c|c|c|}
\hline
\vertical & $x_2$ & $x_{1,1}$ & $x_{1,\unbar}$ & $x_\deuxbar$ & $x_{\unbar,\unbar}$ \\
\hline
\hline
\vertical $\pi_2$               & 1 & . & . & . & . \\
\vertical $\pi_{1,1}$           & 1 & 2 & . & . & . \\
\vertical $\pi_{1,\unbar}$      & 1 & 2 & 2 & . & . \\
\vertical $\pi_\deuxbar$        & 1 & . & . & 2 & . \\
\vertical $\pi_{\unbar,\unbar}$ & 1 & 2 & 4 & 4 & 8 \\
\hline
\end{tabular}\\
\\
\vertical {\bf Table V. Character table of ${\boldsymbol{\Sigma'(W_2)}}$} \\
\end{tabular}}
\end{centerline}

\bigskip

\noindent{\it Remark.} Using these tables, one can check that
$\th_2(x_\deuxbar)(x_{1,1}) = 6 \not= 4 = \th_2(x_{1,1})(x_\deuxbar)$.
In other words, the symmetry property (see \cite{bhs}) does not hold 
in our case.

\bigskip

\begin{eqnarray*}
E_2&=&x_{2}- \frac{1}{2}x_\deuxbar - \frac{1}{4}x_{1,\unbar} +
\frac{1}{4}x_{\bar 1, 1}- \frac{1}{2}x_{1,1}+\frac{1}{4}x_{\bar 1,\bar 1} \\
E_{1,1}&=&\frac{1}{2}\left(x_{1,1} - \frac{1}{2}x_{1,\bar 1} - \frac{1}{2}x_{\bar 1, 1} +
\frac{1}{4}x_{\bar 1,\bar 1} \right)\\
E_{1,\unbar}&=& \frac{1}{2}\left(x_{1,\bar 1} - \frac{1}{2}x_{\bar 1,\bar 1} \right)\\
E_\deuxbar&=& \frac{1}{2}\left(x_{\bar 2} - \frac{1}{2}x_{\bar 1,\bar 1} \right)\\
E_{\unbar,\unbar}&=&\frac{1}{8} x_{\bar 1,\bar 1}\\
\end{eqnarray*}

\begin{centerline}{\bf Table VI. A complete set of
orthogonal primitive idempotents}\end{centerline}

\bigskip

We will now give the Cartan matrix of $\Sigma'(W_2)$.
If $\m \in \Bip(2)$, we denote by $\Pi_\m$ the character of the projective
cover $\QM\Sigma'(W_2) E_\m$ of $\QM_\m$. Write
$$\Pi_\m=\sum_{\mu \in \Bip(2)} \g_{\l\m} \pi_\l.$$
Then $(\g_{\l\m})_{\l,\m \in \Bip(2)}$ is the Cartan matrix of
$\Sigma'(W_2)$. It is given in the following table:

\medskip

\begin{centerline}
{\begin{tabular}{c}
\begin{tabular}{|c|ccccc|}
\hline
\vertical $\lamh~\backslash~\hat{\mu}$ & $(2)$ & $(1,1)$ & $(1,\unbar)$ & $(\deuxbar)$ & $(\unbar,\unbar)$ \\
\hline
\vertical $(2)$                        & 1 & . & . & . & . \\
\vertical $(1,1)$                      & . & 1 & . & . & . \\
\vertical $(1,\unbar)$                 & . & . & 1 & 1 & . \\
\vertical $(\deuxbar)$                 & . & . & . & 1 & . \\
\vertical $(\unbar,\unbar)$            & . & . & . & . & 1 \\
\hline
\end{tabular}\\
\\
\vertical {\bf Table VII. Cartan matrix of ${\boldsymbol{\Sigma'(W_2)}}$} \\
\end{tabular}}
\end{centerline}

\bigskip
\def\matrice#1{\left(\begin{array}{cc} #1 \end{array}\right)}

Let $E_0=E_{1,\unbar}+E_\deuxbar$. Then $(E_2,E_{1,1},E_0,E_\deuxbar)$ is
a complete set of central primitive idempotents (they are of course
orthogonal). Therefore, write $A_\omega = \QM\Sigma'(W_2) E_\omega$,
for $\omega \in \{2, (1,1), 0, \deuxbar\}$. Then
$$\QM\Sigma'(W_2) = A_2 \oplus A_{1,1} \oplus A_\deuxbar \oplus A_0,$$
as a sum of algebras. Morever, $A_2 \simeq \QM$, $A_{1,1} \simeq \QM$,
$A_\deuxbar \simeq \QM$. On the other hand,
$$A_0 = \QM E_{1,\unbar} \oplus \QM E_\deuxbar \oplus
\QM (x_{1,\unbar}-x_{\unbar,1}),$$
as a vector space. Now, let $B$ be the algebra
$$B=\{\matrice{a & b \\ 0 & c\\}~|~a,b,c \in \QM\}.$$
Then the $\QM$-linear map $\s : A_0 \to B$ such that
$$\s(E_{1,\unbar})=\matrice{1 & 0 \\ 0 & 0},\quad
\s(E_\deuxbar)=\matrice{0 & 0 \\ 0 & 1}\quad\text{and}\quad
\s(x_{1,\unbar}-x_{\unbar,1})=\matrice{0 & 1 \\ 0 & 0}$$
is an isomorphism of algebras. Therefore, we have
an isomorphism of algebras
$$\QM\Sigma'(W_2) \simeq \QM \oplus \QM \oplus \QM \oplus B.$$

\bigskip

%
%
%   Biblio
%
%
%

\newpage

%%%%%
%%%%%
%
%
%
%
%
%
% APPENDICE Hohlweg-Baumann
%
%
%
%
%
%
%
%
%%%%%
%%%%%

\section*{Appendix: Comparison with Specht's construction}

\vskip 5mm

\begin{center}\textsc{\large
 Pierre Baumann and Christophe Hohlweg}
\end{center}

\vskip 1cm
The present text is an appendix to the article {\it Generalized  descent
algebra and construction of irreducible characters of  hyperoctahedral
groups}, by C\'edric Bonnaf\'e and the second present author. Our aim
here is to relate two constructions of the irreducible characters of the
hyperoctahedral groups: the one given in that article, and Specht's one
\cite{specht}. Meant as a sequel to Bonnaf\'e and Hohlweg's article, this
text uses the same notations and references.

We first recall briefly Specht's construction, using Macdonald's book
as a reference \cite[I, Appendix B]{macdonald}.

\medskip
\subsection*{Specht's construction.}
Let $G$ be a finite group, let $G_*$ be the set of conjugacy classes
in $G$ and let $G^*$ be the set of irreducible characters of $G$.
Given a conjugacy class $c\in G_*$, we denote by $\zeta_c$ the order
of the centralizer of an element of $c$. We denote the value of a character $\gamma$ of $G$ at any element of a comjugacy class $c\in G_*$ by $\gamma(c)$.

We denote the wreath product $G\wr\mathfrak S_n$ by $G_n$. This wreath
product is the semidirect product $G^n\rtimes\mathfrak S_n$ for the
action of $\mathfrak S_n$ on $G^n$ given by
$$\sigma\cdot(g_1,\ldots,g_n)=(g_{\sigma^{-1}(1)},\ldots,
g_{\sigma^{-1}(n)}),$$
where $\sigma\in\mathfrak S_n$ and $(g_1,\ldots,g_n)\in G^n$, so that we
can always represent an element in $G_n$ as a product
$(g_1\ldots,g_n)\;\sigma$.

Given a complex representation $V$ of $G$, we construct a complex
representation $\eta_n(V)$ of $G_n$ on the space $V^{\otimes n}$ by
letting a product $(g_1\ldots,g_n)\;\sigma$ acting on a pure tensor
$v_1\otimes\cdots\otimes v_n\in V^{\otimes n}$ in the following way:
$$((g_1\ldots,g_n)\;\sigma)\cdot(v_1\otimes\cdots\otimes v_n)=(g_1\cdot
v_{\sigma^{-1}(1)})\otimes\cdots\otimes(g_n\cdot v_{\sigma^{-1}(n)}).$$
The character of $\eta_n(V)$ does not depend on $V$ but only of its
character; if $\rho$ denotes the latter, then we will denote the former
by $\eta_n(\rho)$.

We let $\mathcal P$ be the set of all partitions, and we set
$\mathcal P_G=\mathcal P^{G^*}$. Given an element $\lambda=
(\lambda_\gamma)_{\gamma\in G^*}$ in $\mathcal P_G$, we denote by
$|\lambda |$ the sum $\sum_{\gamma}|\lambda_\gamma|$.

Now let $\Lambda_{\mathbb C}$ be the (free) ring of symmetric polynomials with
complex coefficients. As is well-known, $\Lambda_{\mathbb C}$ is
generated over $\mathbb C$ by a countable family of algebraically
independant elements: one can choose for generators the family
$(h_n)_{n\geq1}$ of complete symmetric functions or the family
$(p_n)_{n\geq1}$ of power sums. On the other hand, the family of
Schur functions $(s_\lambda)_{\lambda_\in\mathcal P}$ is a basis of
the vector space $\Lambda_{\mathbb C}$. Following Macdonald, we denote
by $\Lambda_{\mathbb C}(G)$ the ring of polynomials over $\mathbb C$ in
the family of variables $\bigl(p_n(c)\bigr)_{n\geq1,c\in G_*}$. Setting
$$p_n(\gamma)=\sum_{c\in G_*}\zeta_c^{-1}\gamma(c)p_r(c)$$
for any $\gamma\in G^*$, one can easily check that $\Lambda_{\mathbb C}(G)$
is also the ring of polynomials in the variables $\bigl(p_n(\gamma)\bigr)_{
n\geq1,\gamma\in G^*}$. Every symmetric polynomial $P\in\Lambda_{\mathbb C}$
can be expressed as a polynomial with complex coefficients in the power
sums $p_n$; given $\gamma\in G^*$, we denote by $P(\gamma)$ the element of
$\Lambda_{\mathbb C}(G)$ obtained by replacing the variables $p_n$ by the
variables $p_n(\gamma)$ in the expression of $P$. Given an element
$\lambda=(\lambda_\gamma)_{\gamma\in G^*}$ in $\mathcal P_G$,
we set
$$s_{\lambda}=\prod_{\gamma\in G^*}s_{\lambda_\gamma}(\gamma).$$

The set of complex irreducible characters of $G_n$ is a basis of the
algebra of complex-valued class functions of $G_n$, so that we can
denote this latter by $\mathbb C\Irr(G_n)$. The direct sum
$$R(G)=\bigoplus_{n\geq0}\mathbb C\Irr(G_n)$$
can then be endowed with the structure of a commutative and cocommutative
$\mathbb N$-graded Hopf algebra, where the product is given by (the maps induced on the level of characters by) the induction
functors $\Ind_{G_m\times G_n}^{G_{m+n}}$ and the coproduct is afforded
likewise by the restriction functors $\Res_{G_m\times G_n}^{G_{m+n}}$
\cite[I, Appendix~B, 4 and I, 7, Example 26]{macdonald}. Since
$\Lambda_{\mathbb C}(G)$ is a free commutative algebra, there is a
unique homomorphism of $\mathbb C$-algebras
$$\ch^{-1}:\Lambda_{\mathbb C}(G)\to R(G)$$
with the following property: for each $n\geq0$ and each $c\in G_*$,
$\ch^{-1}$ maps the variable $p_n(c)$ to the characteristic function
of the conjugacy class of $G_n$ consisting of the products $(g_1,\ldots,g_n)
\;\sigma$, where the permutation $\sigma\in\mathfrak S_n$ is a $n$-cycle
and the product $g_1g_2\cdots g_n$ belongs to the conjugacy class $c$.
It turns out that $\ch^{-1}$ is an isomorphism of Hopf algebra, whose
inverse will be denoted by $\ch$. Then, using arguments of orthogonality
and integrality, it can be shown \cite[I, Appendix B, 9]{macdonald} that
the image under $\ch$ of the irreducible characters of $G_n$ are the
elements $s_\lambda$, where $\lambda\in\mathcal
P_G$ is such that $|\lambda |=n$.

Later on, we will need to know the image under $\ch$ of characters
$\eta_n(\rho)$. We do the computation now.
\begin{lem}\label{le:CalcChEta}
Let $\gamma_1$, \dots, $\gamma_s$ the irreducible characters of $G$, let
$c_1$, \dots, $c_s$ be non-negative integers, and set $\rho=c_1\gamma_1+
\cdots+c_s\gamma_s$. Then
$$\sum_{n\geq0}\ch\bigl(\eta_n(\rho)\bigr)=\prod_{i=1}^s\Biggl(\sum_{n\geq0}
h_n(\gamma_i)\Biggr)\biggr.^{c_i}.$$
\end{lem}

\begin{proof}
The proof given in \cite[I, Appendix B, 8]{macdonald} for the case where
$\rho$ is irreducible can be easily adapted. Indeed in the computation
that follows Equation (8.2) in that reference, the steps that lead to the
equality
$$\sum_{n\geq0}\ch\bigl(\eta_n(\gamma)\bigr)=\exp\Biggl(\sum_{r\geq1}
\frac1r\sum_{c\in G_*}\zeta_c^{-1}\gamma(c)p_r(c)\Biggr)$$
are valid even if the character $\gamma$ is reducible. Applying this
formula to the character $\rho$, we get
\begin{align*}
\sum_{n\geq0}\ch\bigl(\eta_n(\rho)\bigr)&=\exp\Biggl(\sum_{r\geq1}\frac1r
\sum_{c\in G_*}\zeta_c^{-1}\rho(c)p_r(c)\Biggr)\\&=\prod_{i=1}^s
\Biggl[\exp\Biggl(\sum_{r\geq1}\frac1r\sum_{c\in G_*}\zeta_c^{-1}
\gamma_i(c)p_r(c)\Biggr)\Biggr]\biggr.^{c_i}\\&=\prod_{i=1}^s
\Biggl[\exp\Biggl(\sum_{r\geq1}\frac1rp_r(\gamma_i)\Biggr)\Biggr]
\biggr.^{c_i}\\&=\prod_{i=1}^s\Bigl[\sum_{n\geq0}h_n(\gamma_i)\Bigr]^{c_i},
\end{align*}
the last step in the computation coming from Newton's formulas.
\end{proof}

\medskip
\subsection*{The comparison result.}
Having now recalled Specht's construction of the irreducible characters
for the wreath product $G\wr\mathfrak S_n$ of an arbitrary finite group
$G$ by the symmetric group $\mathfrak S_n$, we can specialize to the
case where $G$ is the group $W=\mathbb Z/2\mathbb Z$ with two elements.
The notation $W_n$ for the wreath product $W\wr\mathfrak S_n$ then agrees
with its use by Bonnaf\'e and Hohlweg. The Hopf algebra $R(W)$ is
identical to the complexified Hopf algebra $\mathcal{CHAR}
\otimes_{\mathbb Z}\mathbb C$. The set $W^*$ of irreducible characters
of $W$ has two elements, namely the trivial character $\tau$ and the
signature $\varepsilon$. One can view an element $\lambda=
(\lambda_\tau,\lambda_\varepsilon)$ of $\mathcal P_W$ as a bipartition
$(\lambda^+,\lambda^-)$ by setting $\lambda^+=\lambda_\tau$ and
$\lambda^-=\lambda_\varepsilon$. As a final piece of notation, we set
$\lambda^*=\bigl(\lambda^+,(\lambda^-)^t\bigr)$ for any
bipartition $\lambda=(\lambda^+,\lambda^-)$.

Generalizing Poirier and Reutenauer's work \cite{poirier-reutenauer} for
symmetric groups to the case of $W_n$, we define a linear map:
\begin{eqnarray*}
f:\,\cop\otimes_{\mathbb Z}\mathbb C&\longrightarrow&\Lambda_{\mathbb C}(W)
\end{eqnarray*}
by setting $f(z_Q)=s_{(\sh Q)^*}$ for any bitableau $Q$.
With all these notations, our result can be stated as follows:

\begin{thm}\label{appendix}
The following diagram of $\mathbb N$-graded
Hopf algebras
$$\xymatrix{\cop\otimes_{\mathbb Z}\mathbb C\ar@{->>}[r]^f
\ar@{->>}[dr]^{\widetilde{\theta}}&\Lambda_{\mathbb C}(W)\\
\Sigma'\otimes_{\mathbb Z}\mathbb C\ar@{^{(}->}[u]^i\ar@{->>}[r]_\theta&
\mathcal{CHAR}\otimes_{\mathbb Z}\mathbb C\ar[u]^\sim_\ch}$$
is commutative. In particular $\ch(\xi_{\lambda})=s_{
\lambda^*}$, for any bipartition $\lambda$, so that Bonnaf\'e
and Hohlweg's construction is equivalent to Specht's one, up to a relabelling.
\end{thm}

\medskip
Some further notation and a bijection will be needed for the proof. We
present them now.

\subsection*{Some notations and a bijection.}
We call quasicomposition a sequence $E=(e_1,e_2,e_3,\ldots)$ of
non-negative integers, all of whose terms but a finite number vanish.
The size $|E|$ of $E$ is the sum $e_1+e_2+e_3+\cdots$ of the terms.
Given a partition $\mu$ and a quasicomposition $E$, we denote by
$\tab\mu E$ the set of all semistandard tableau of shape $\mu$ and
weight $E$, that is the set of all fillings of the Ferrers diagram of
shape $\mu$ with positive integers, in such a way that the numbers are
weakly increasing from left to right in the rows, strictly increasing
from top to bottom in the columns, and that there is $e_1$ times the
number $1$, $e_2$ times the number $2$, and so on \cite[p.~5]{macdonald}.
The set $\tab\mu E$ is of course empty unless $|\mu|=|E|$. Given any
quasicomposition $E=(e_1,e_2,e_3,\ldots)$, the formula
$$h_{e_1}\,h_{e_2}\,h_{e_3}\cdots=\sum_{\mu\in\mathcal P}
\bigl|\tab\mu E\bigr|s_\mu$$
holds in $\Lambda_{\mathbb C}$ (see \cite[I, (6.4)]{macdonald} for a proof).

Now we fix a positive integer $n$ and a signed composition
$C=(c_1,\ldots,c_\ell)$ of it. Let $\ell$ be the length of $C$.
We define $\Comp(C)$ as the set of all quasicompositions
$D=(d_1,\ldots,d_\ell)$ such that $d_i=0$ if $c_i>0$ and $0\leq d_i\leq
-c_i$ if $c_i<0$. Given such a $D$, we further define two quasicompositions
$T_{C,D}=(t_1,\ldots,t_\ell)$ and $E_{C,D}=(e_1,\ldots,e_\ell)$ by
$$t_i=\begin{cases}c_i&\text{if $c_i>0$,}\\d_i&\text{if $c_i<0$;}
\end{cases}\quad\text{and}\quad e_i=\begin{cases}0&\text{if $c_i>0$,}\\
-c_i-d_i&\text{if $c_i<0$.}
\end{cases}$$
The signed composition obtained by omitting the zeros in the list
$$(-e_1,t_1,-e_2,t_2,\ldots,-e_\ell,t_\ell)$$
will be denoted by $B_{C,D}$. For instance, for
$C=(2,\bar2,\bar3,1,\bar1,2,2,\bar2)\VDash15$, we can choose
$D=(0,0,2,0,1,0,0,0)$, and then $T_{C,D}=(2,0,2,1,1,2,2,0)$,
$E_{C,D}=(0,2,1,0,0,0,0,2)$ and $B_{C,D}=(2,\bar2,\bar1,2,1,1,2,2,\bar2)$.

Finally, given a bipartition $\lambda=(\lambda^+,\lambda^-)$
and a signed composition $C$ with $|\lambda|=|C|$, we define
$\bitab{\lambda}C$ as the set of all standard bitableaux $Q$
such that $\sh(Q)=\lambda^*$ and $C\leftarrow\mathbf C(Q)$ (see Remark~\ref{composition tableau}).

\smallskip
One of the key to the proof of Theorem~\ref{appendix} is the following
combinatorial result.
\begin{prop}\label{bijection annexe}
Given a bipartition $\lambda$ and a signed composition $C$ with
$|\lambda |=|C|$, the sets $\bitab{\lambda}C$ and
$$\coprod_{D\in\Comp(C)}\tab{\lambda^+}{T_{C,D}}\times
\tab{\lambda^-}{E_{C,D}}$$
have the same cardinality.
\end{prop}

\begin{proof}
Let $n$ be a positive integer, $C$ be a signed composition of $n$,
and $\lambda=(\lambda^+,\lambda^-)$ be a bipartition with
$|\lambda|=n$. We construct inverse bijections between
$\bitab{\lambda}C$ and
$$\coprod_{D\in\Comp(C)}\tab{\lambda^+}{T_{C,D}}\times
\tab{\lambda^-}{E_{C,D}}$$
as follows.

First let $(R,S)$ be in the second set, so that
$R\in\tab{\lambda^+}{T_{C,D}}$ and $S\in\tab{\lambda^-}{E_{C,D}}$
for some $D\in\Comp(C)$.  We can put a total order on the boxes in $R$
and $S$ by requiring that:
\begin{itemize}
\item A box is smaller than another one if the label written in it is
smaller than the one in the other.
\item Given two boxes with the same label in it, a box in $S$ is
smaller than a box in $R$.
\item For boxes containing the same label and located in the same
tableau ($R$ or $S$), boxes located south-west are smaller than boxes
located north-east.
\end{itemize}
We can then enumerate in increasing order the boxes in $R$ and $S$.
Filling now each box of $R$ and $S$ by its rank of appearance in the
enumeration, we construct a standard bitableau $\tilde Q$ of shape
$\lambda$. We then define $Q$ as the bitableau obtained from
$\tilde Q$ by transposing $\tilde Q^-$, so that $Q$ has shape
$\lambda^*$. Comparing this construction with the combinatorial rule
in Remark~\ref{composition tableau} that computes $\bs C (Q)$, we easily check that
the signed composition $B_{C,D}$ can be obtained from $\mathbf C(Q)$ by
refinement of the parts, so that $C{\buildrel B\over\longleftarrow}
B_{C,D}{\buildrel R\over\longleftarrow}\mathbf C(Q)$, which implies
$Q\in\bitab{\lambda}C$.

In the other direction, let $Q$ be a given element in
$\bitab{\lambda}C$. From Theorem~\ref{refinement signed compositions}, there
exists a unique signed composition $B$ such that $C{\buildrel B\over
\longleftarrow}B{\buildrel R\over\longleftarrow}\mathbf C(Q)$, and we can
find a (unique) element $D\in\Comp(C)$ so that $B=B_{C,D}$. Now we
transpose $Q^-$ and get a bitableau $\tilde Q$. We construct a list
$L=(l_1,l_2,\ldots,l_n)$ of positive integers by placing first $|c_1|$
times the number $1$, then $|c_2|$ times the number $2$, and so on. Then
we substitute $l_1$ to $1$, $l_2$ to $2$, and so on, in the boxes of
$\tilde Q$, and obtain in this way a pair $(R,S)$ of tableaux of shapes
$\lambda^+$ and $\lambda^-$ respectively. The fact that $B_{C,D}
{\buildrel R\over\longleftarrow}\mathbf C(Q)$ implies that this
construction yield two semistandard tableaux $R$ and $S$ with weights
$T_{C,D}$ and $E_{C,D}$ respectively, that is to say
$$(R,S)\in\tab{\lambda^+}{T_{C,D}}\times\tab{\lambda^-}{E_{C,D}}.$$

It is a routine task to check that the two above constructions yield
mutually inverse bijections between $\coprod_{D\in\Comp(C)}
\tab{\lambda^+}{T_{C,D}}\times\tab{\lambda^-}{E_{C,D}}$ and
$\bitab{\lambda}C$.
\end{proof}

We end this paragraph by an example that illustrates the constructions
needed in the proof above. We take $n=15$ and choose the same signed
composition $C$ as in the previous example, namely
$$C=(2,\bar2,\bar3,1,\bar1,2,2,\bar2).$$
We choose $\lambda^+=631$ and $\lambda^-=41$, so that
$\lambda^*=(631,21^3)$. Starting from the pair $(R,S)$ with
$$R=\young(113347,567,6)\qquad\text{and}\qquad S=\young(2238,8),$$
we construct $\tilde Q=(\tilde Q^+,\tilde Q^-)$ where
$$\tilde Q^+=\young(12678\thirteen,9\eleven\twelve,\ten)\qquad\text{and}
\qquad\tilde Q^-=\young(345\fifteen,\fourteen),$$
whence $Q=(Q^+,Q^-)$ with
$$Q^+=\tilde Q^+=\young(12678\thirteen,9\eleven\twelve,\ten)\qquad
\text{and}\qquad Q^-={}^t\tilde Q^-=\young(3\fourteen,4,5,\fifteen).$$
Since $\mathbf C(Q)=(2,\bar3,3,1,4,\bar2)$, it holds that
$C{\buildrel B\over\longleftarrow}B{\buildrel R\over\longleftarrow}
\mathbf C(Q)$ with
$$B=(2,\bar2,\bar1,2,1,1,2,2,\bar2),$$
which implies $C\leftarrow\mathbf C(Q)$.

In the other direction, we start from the bitableau $Q$. We observe
that the signed composition $B$ such that $C{\buildrel B\over
\longleftarrow}B{\buildrel R\over\longleftarrow}\mathbf C(Q)$
is $B_{C,D}$, where $D$ is given by $D=(0,0,2,0,1,0,0,0)$. Now we write
down the list
$$L=(1,1,2,2,3,3,3,4,5,6,6,7,7,8,8)$$
from $C$. Transposing the negative tableau $Q^-$, we write down $\tilde Q$
and substitute the elements of $L$ to the numbers in the boxes of
$\tilde Q$. We recover our original pair $(R,S)$. We easily verify that
$R$ has weight
$$T_{C,D}=(2,0,2,1,1,2,2,0)$$
and that $S$ has weight
$$E_{C,D}=(0,2,1,0,0,0,0,2).$$

\subsection*{Proof of Theorem~\ref{appendix}.}\

\noindent{\bf1}\quad We first compute the image by $\ch$ of the induced
character $\Ind_{\mathfrak S_n}^{W_n}1$ of $W_n$, where $n$ is a positive
integer. To do that, we construct the complex representation $\eta_n(V)$
of $W_n$, where $V$ is the left regular representation of $W=\mathbb Z/2\mathbb Z$.
Denoting by $\mathbb C_1$ the trivial representation of $\mathfrak S_n$,
we then observe that the isomorphism of vector spaces $\Ind_{\mathfrak
S_n}^{W_n}\mathbb C_1\cong\eta_n(V)$ given by the sequence of natural
identifications
$$\Ind_{\mathfrak S_n}^{W_n}\mathbb C_1\cong\mathbb CW_n\otimes_{\mathbb
C\mathfrak S_n}\mathbb C_1\cong\mathbb C(W^n)\cong(\mathbb CW)^{\otimes n}
=V^{\otimes n}=\eta_n(V)$$
is compatible with the action of $W_n$. Since $V$ has $\tau+\varepsilon$
for character, it follows that $\Ind_{\mathfrak S_n}^{W_n}1=\eta_n
(\tau+\varepsilon)$. Lemma~\ref{le:CalcChEta} now implies that
$$\sum_{n\geq0}\ch\bigl(\Ind_{\mathfrak S_n}^{W_n}1\bigr)=\Biggl(
\sum_{n\geq0}h_n(\tau)\Biggr)\Biggl(\sum_{n\geq0}h_n(\varepsilon)\Biggr).$$

On the other side, it is easy to check that $\eta_n(\tau)$ is the trivial
character of $W_n$. Therefore $\ch$ maps the trivial character $\Ind_{W_n}^{W_n}1$ of $W_n$  to $h_n(\tau)$. To comply
with the philosophy used by Bonnaf\'e and Hohlweg, we will write
for any positive integer $n$
$$\varphi_{\pm n}=\ch\bigl(\Ind_{W_{\pm n}}^{W_n}1\bigr)=\begin{cases}
\ch\bigl(\Ind_{W_n}^{W_n}1\bigr)=h_n(\tau)&\text{for '$+$' sign,}\\[5pt]
\ch\bigl(\Ind_{\mathfrak S_n}^{W_n}1\bigr)=\sum_{k=0}^n
h_k(\tau)h_{n-k}(\varepsilon)&\text{for '$-$' sign.}\end{cases}$$

\medskip
\noindent{\bf2}\quad We now prove the equality $f\circ i=\ch\circ\theta$.
Given any signed composition $C=(c_1,\ldots,c_\ell)$, there holds
$x_C=x_{c_1}\cdots x_{c_\ell}$. Since $\theta$ is a morphism of Hopf
algebras, we can write
$$\Ind_{W_C}^{W_{|C|}}1_C =\theta(x_C)=\theta(x_{c_1})\cdots\theta(x_{c_\ell})
=\Ind_{W_{c_1}}^{W_{|c_1|}}1\CLprod\cdots\CLprod
\Ind_{W_{c_\ell}}^{W_{|c_\ell|}}1,$$
and taking its image under $\ch$,
$$\ch\Ind_{W_C}^{W_{|C|}}=\ch\circ\theta(x_C)=\varphi_{c_1}\cdots
\varphi_{c_\ell}.$$
The formula
$$\varphi_{-n}=\sum_{k=0}^nh_k(\tau)h_{n-k}(\varepsilon),$$
valid for any positive integer $n$, makes possible to continue the
computation:
$$\ch\circ\theta(x_C)=\sum_{D\in\Comp(C)}h_{t_1}(\tau)\cdots h_{t_\ell}
(\tau)\;h_{e_1}(\varepsilon)\cdots h_{e_\ell}(\varepsilon),$$
where the quasicompositions $(t_1,\ldots,t_\ell)$ and $(e_1,\ldots,
e_\ell)$ appearing in the sum are $T_{C,D}$ and $E_{C,D}$ respectively.
We thus get, using Proposition~\ref{bijection annexe} and the decomposition of $X_C$ given at the end of Remark~\ref{composition tableau}:
\begin{align*}
\ch\circ\theta(x_C)&=\sum_{D\in\Comp(C)}\left(\sum_{\lambda^+\in\mathcal P}
\Bigl|\tab{\lambda^+}{T_{C,D}}\Bigr|s_{\lambda^+}(\tau)\right)\;
\left(\sum_{\lambda^-\in\mathcal P}\Bigl|\tab{\lambda^-}
{E_{C,D}}\Bigr|s_{\lambda^-}(\varepsilon)\right)\\
&=\sum_{(\lambda^+,\lambda^-)\in\mathcal P_W}\left(\sum_{D\in\Comp(C)}
\Bigl|\tab{\lambda^+}{T_{C,D}}\times\tab{\lambda^-}{E_{C,D}}\Bigr|\right)\;
s_{\lambda^+}(\tau)\;s_{\lambda^-}(\varepsilon)\\
&=\sum_{\lambda\in\mathcal P_W}\Bigl|\bitab{\lambda}
C\Bigr|s_{\lambda}\\
&=\sum_{\substack{Q\text{ std. bitableau}\\C\gets\bs C(Q)}}
s_{(\sh Q)^*}\\
&=\sum_{\substack{Q\text{ std. bitableau}\\C\gets\bs C(Q)}}
f(z_Q)\\
&=f\circ i(x_C).
\end{align*}
Since the elements $x_C$ generate $\Sigma'\otimes_{\mathbb Z}\mathbb C$ as
a vector space, it follows that $\ch\circ\theta=f\circ i$.

\medskip
\noindent{\bf3}\quad To complete the proof, it now suffices to show that
$f=\ch\circ\tilde\theta$. We first observe that both members of this
equality coincide on the image of $i$ in $\cop\otimes_{\mathbb Z}
\mathbb C$, since
$$\ch\circ\theta=f\circ i\qquad\text{and}\qquad\theta=\tilde\theta
\circ i.$$
On the other hand, $f$ and $\ch\circ\tilde\theta$ have the same kernel,
namely the vector space $\QC_n^\perp$ spanned by the elements $z_Q-z_{Q'}$ for
standard bitableaux $Q$ and $Q'$ of the same shape (see Theorem~\ref{theoreme theta}). Since $\theta$
is surjective, this kernel, together with the image of $i$, spans
$\cop\otimes_{\mathbb Z}\mathbb C$. The result follows easily from
these facts.
\end{document}